\documentclass[12pt,a4paper]{article}
\usepackage{latexsym}
\usepackage{amsmath}
\usepackage{amssymb}
\usepackage{amsfonts}
\usepackage{verbatim}
\usepackage[latin1]{inputenc}
\newtheorem{Th}{Theorem}[section]

\newtheorem{proposition}[Th]{Proposition}
\newtheorem{lemme}[Th]{Lemma}
\renewcommand{\theequation}{\arabic{section}.\arabic{equation}}
\let \ssection=\section
\renewcommand{\section}{\setcounter{equation}{0}\ssection}

\def\^#1{\if#1i{\accent"5E\i}\else{\accent"5E #1}\fi}
\def\"#1{\if#1i{\accent"7F\i}\else{\accent"7F #1}\fi}
\allowdisplaybreaks[4]

\newcommand{\E}{\mathbb E}

\newcommand{\NN}{\mathbb N}
\newcommand{\RR}{\mathbb R}

\newcommand{\CC}{\mathcal C}
\newcommand{\HH}{\mathcal H}
\newcommand{\FF}{\mathcal F}

\newcommand{\ph}{\varphi}
\newcommand{\dd}{\partial}

\newcommand{\ddd}{\displaystyle}
\newcommand{\dint}{\int\hspace{-2mm}\int} 
\newcommand{\kk}{\underline{k}}
\newcommand{\xx}{\underline{x}}

\begin{document}
\title{On implicit and explicit discretization schemes for parabolic SPDEs in any dimension}
\author{Annie MILLET$^*$ \& Pierre-Luc MORIEN$^{**}$}

\date{}
\def\^#1{\if#1i{\accent"5E\i}\else{\accent"5E #1}\fi}
\def\"#1{\if#1i{\accent"7F\i}\else{\accent"7F #1}\fi}
\maketitle
 \begin{center}
{\it  * Laboratoire de Probabilit\'es, Universit\'e Paris 6, 4 place Jussieu\\
 75252 Paris, France\; and
 \\ SAMOS-MATISSE, Universit\'e Paris 1, 90 Rue de Tolbiac, 75634 Paris Cedex 13 France
 \\

 ** MODAL'X,  Universit\'e Paris 10, 200 avenue de la R\'epublique \\
 92001 Nanterre Cedex, France
\\

 E-mail addresses: \quad * amil@ccr.jussieu.fr \quad ** morien@u-paris10.fr}
  \end{center}

\begin{abstract}
We study the speed of convergence of the explicit and implicit
space-time discretization schemes of the solution $u(t,x)$ to a
parabolic partial differential equation in any dimension perturbed
by a space-correlated Gaussian noise. The coefficients only depend
on $u(t,x)$ and the influence of the correlation on the speed is
observed.  
\smallskip

\smallskip

\end{abstract}

\section{Introduction}\label{intro}
Discretization schemes for parabolic SPDEs driven by the space-time white
noise have been considered  by several authors. I.~Gy\"ongy and D.~Nualart \cite{GN1} and \cite{GN2},
 have studied implicit time discretization schemes for
the heat equation in dimension 1. J.~Printems \cite{P1} has studied several time discretization
schemes (implicit and  explicit Euler schemes as well as the Crank-Nicholson one)
for Hilbert-valued parabolic SPDEs, such
as the Burgers equation on $[0,1]$, introduced several notions of order of convergence
in order to deal with coefficients with polynomial growth and  proved  convergence in
the Hilbert space norm. This work  has been completed
by E. Hausenblas \cite{H}, who studied several schemes for quasi-linear equations
driven by  a nuclear noise, and taking values in a Hilbert  or a Banach space $X$.
Several approximation procedures (such as the Galerkin approximation, finite difference methods
or  wavelets approximations)
were considered, but the coefficients of the SPDE were  supposed to
depend on the whole function $u(t,.)$  in $X$, and not only on 
 $(t,x)$. Notice that, unlike \cite{H},  the coefficients considered  in this paper do not
depend on the whole function $u(s,.)$.

I.~Gy\"{o}ngy  \cite{G1} has studied the strong speed of convergence in the norm of uniform convergence
over the  space variable
for a space finite-difference scheme $u^n$ with mesh $1/n$ for the parabolic SPDE with homogeneous Dirichlet's
boundary conditions. He has also studied the speed of convergence of
an implicit (resp. explicit) finite-difference discretization scheme $u^{n,m}$ (resp. $u_m^n$) with time mesh
$T/m$ and space mesh $1/n$ for the solution $u$ to the following parabolic SPDE in dimension 1
driven by the space-time white noise $W$:
\begin{equation}\label{parablanc}
\left\{\begin{array}{l}
\frac{\partial u}{\partial t}(t,x)=\frac{\partial^2u}{\partial x^2}(t,x)+
\sigma(t,x,u(t,x))\frac{\partial^2 W}{\partial t\partial x}+b(t,x,u(t,x))\, ,\\
u(t,0)=u(t,1)=0 \, ,
\end{array}\right.
\end{equation}
with the initial condition $u_0$.  He has proved that, if the coefficients $\sigma(t,x,.)$
and $b(t,x,.)$ satisfy the usual Lipschitz property uniformly in $(t,x)$
and if the functions
$\sigma(t,x,y)$ and $b(t,x,y)$ are $1/4$-H\"{o}lder continuous in $t$ and $1/2$-H\"{o}lder
continuous in $x$ uniformly with respect to the other variables, then for $t\in ]0,T]$,
$p\in [1,+\infty[$, $0<\beta<\frac{1}{4}$ and $0<\gamma<\frac{1}{2}$  one has:
\begin{equation}\label{vitesblanc}
\sup_{x\in [0,1]} \E\left( \left| u^{n,m}(t,x)-u(t,x)\right|^p\right)\leq K(t)
\left( m^{-\beta p}+n^{-\gamma p}\right)\, .
\end{equation}
Furthermore,
if $u_0\in \CC^3([0,1])$,
then (\ref{vitesblanc}) holds on $[0,T]$ with $\beta=\frac{1}{4}$,  $\gamma=\frac{1}{2}$
and  with a constant $K$ which does not depend on $t$. A similar
result holds for the explicit scheme $u_m^n$ if $\frac{n^2 T}{m}\leq q<\frac{1}{2}$.

A. Debussche and J. Printems \cite{P2} have implemented simulations of a discretization scheme for the KDV equation,
and C. Cardon-Weber \cite{CW} has studied explicit and implicit
discretization schemes for the function-valued solution to the stochastic 
Cahn-Hilliard equation in dimension $d\leq 3$ when the driving noise
 is the space-time white  noise. 
The polynomial growth  of the drift term made her  require  the
diffusion coefficient $\sigma$  to be  bounded, and
she  proved convergence in probability  (respectively in $L^p$
 with a given rate of a localized version) of the scheme.

In the present  paper, we deal with a $d$-dimensional version of (\ref{parablanc}).
As it is well-known, we can no longer use the space-time white noise for the
perturbation; indeed, in dimension $d\geq 2$, the Green function associated with
$\frac{\partial}{\partial t}-\Delta$ with the homogeneous Dirichlet
boundary conditions on $[0,1]^d$ is not square integrable. Thus, we replace $W$ by some
Gaussian process $F$  which is white in time and has  a space correlation given by a Riesz
potential $f(r)=r^{-\alpha}$, i.e., such that if 
$A$ and $B$ are bounded Borel subsets of ${\mathbb R}^d$, 
$ E\big( F(s,A)\, F(t,B)\big)= (s\wedge t)\,\int_A dx \int_B dy |x-y|^{-\alpha} $  for some $\alpha \in ]0,2\wedge d[$.
See e.g. \cite{KZ}, \cite{D}, \cite{PZ} and \cite{CWM} for more general
results concerning necessary and sufficient conditions on the covariance of
the Gaussian noise $F$ ensuring the existence of a function-valued solution
to (\ref{parablanc}) with $F$ instead of $W$.

The aim of this paper is threefold. We at first study the speed of convergence
of  space and space-time finite discretization implicit (resp. explicit)
schemes in dimension  $d\geq 1$,  i.e.,  on the grid
 $(\frac{iT}{m},(\frac{j_k}{n}, 1\leq k\leq d))$, $0\leq i\leq m$, $0\leq j_k\leq n$ and extended to
$[0,T]\times [0,1]^d$ by linear interpolation. As in \cite{G1} and \cite{G2}, the processes
$u^n$ and $u^{n,m}$ (resp. $u_m^n$) have an evolution formulation written in terms of
 approximations $(G_d)^n$,
$(G_d)^{n,m}$ and $(G_d)_m^n$  of the Green function $G_d$, while $u$ is solution of an evolution
 equation defined in terms of $G_d$. These  evolution equations
involve   stochastic integrals
 with respect to  the worthy martingale-measure defined by $F$ (see e.g.
\cite{W} and \cite{D}).
As usual, the speed of convergence is given by the norm of the
differences of stochastic integrals; more
precisely, the optimal speed of convergence for the implicit scheme is  the norm of the
difference $G_d(.,x,.)-(G_d)^{n,m}(.,x,.)$ in $L^2([0,T], \HH_d)$, where $\HH_d$ is
the Reproducing Kernel  Hilbert Space defined by the covariance function.
More precisely, if $\varphi$ and $\psi$ are
continuous functions on $Q=[0,1]^d$, set
\begin{equation}\label{RKHS}
<\varphi,\psi>_{\HH_d} = \int_Q\int_Q \varphi(x)\, f(|x-y|)\, \psi(y) dx\, dy\, .
\end{equation}
We denote by $\HH_d$  the completion of this pre-Hilbert space; note that $\HH_d$
elements which are not functions and that a function $\varphi $ belongs to
$\HH_d$ if and only if $\int_Q\int_Q |\varphi(y)|\, f(|y-z|)\, |\varphi(z)|\, dy\, dz<+\infty$.
However, unlike in  \cite{G1} and \cite{G2}, the functions
$ 
 \varphi_j(x) = \sqrt{2}\,
\sin(j\pi x)$, $j\geq 1$  and 
$\varphi_j(\kappa_n(x))$, $ 1\leq j\leq n$, 
where $\kappa_n(y)=[ny]\, n^{-1}$  are not an orthonormal family of $\HH_1$. Thus, even in dimension $d=1$, the use of the
Parseval identity has to be replaced by more technical computations based on Abel's summation method.
Similar results could be obtained for a more general space covariance, provided that it is
absolutely continuous and that its density $f$ satisfies some integrability property at the origin
(see e.g. \cite{D}, \cite{PZ}). However, the
speed of convergence would depend on integrals including $f$, which would make  the results  less transparent
 than that stated in the case of Riesz potentials.
The key technical lemmas, giving upper estimates of
$ \| G_d(.,x,.)-(G_d)^n(.,x,.)\|_{L^2([0,\infty[,\HH_d)}$
and $\| (G_d)^n(.,x,.)-(G_d)^{n,m}(.,x,.)\|_{L^2([0,T],\HH_d)}$ (resp.
$\| (G_d)^n(.,x,.)-(G_d)^n_m(.,x,.)\|_{L^2([0,T],\HH_d)}$), are proved in section \ref{cruciale}.

We describe the discretization schemes in any dimension $d\geq 1$
and introduce some notations in section 2.
 In section \ref{convergence},  an argument similar to that in \cite{G1} 
 shows that for $0< \alpha <d\wedge 2$,
 and $p\in[1,+\infty[$, if $u_0$ is regular enough, then
\begin{equation}\label{vitesn1}
\sup_{(t,x)\in [0,+\infty[\times Q}\;
\E(\| u(t,x)-u^n(t,x)\|^{2p})\leq C_{p,\alpha} \, n^{-(2-\alpha)p}\,, 
\end{equation}
and extending \cite{G2} we prove in section 4 that
\begin{equation}\label{vitesm1}
\sup_{(t,x)\in [0,T]\times Q}\; \E(\| u^n(t,x)-u^{n,m}(t,x)\|^{2p})\leq C_{p,\alpha} \, m^{-(1-\frac{\alpha}{2})p}\, .
\end{equation}
If $d=1$, as  $\alpha \nearrow 1$ the space density becomes more and more degenerate and
the speed of convergence approaches  that obtained by Gy\"ongy for the  space-time white 
noise. 

In dimension $d\geq 2$, the proof depends on the product form
of the Green function and its approximations, as well as of upper estimates
of $|x-y|^{-\alpha}$
in terms of $\prod_{i=1}^d |x_i-y_i|^{-\alpha_i}$ for some well-chosen $\alpha_i$.
 Thus, estimates of the $\HH_d$-norm of
the differences of $G_d(s,x,.)-(G_d)^n(s,x,.)$, $(G_d)^n(s,x,.)-(G_d)^{n,m}(s,x,.)$
and $(G_d)^n(s,x,.)-(G_d)_m^n(s,x,.)$ in dimension $d\geq 2$
depend  on bounds of the $\HH_1$-norm of similar  differences as well as of
$\HH_r$-norms of $G(s,x,.)$, $G^n(s,x,.)$ and $G^{n,m}(s,x,.)$ for $r<d$.

Section \ref{numerique} contains some numerical results. For $T=1$,
we have implemented in C the
(more stable) implicit discretization scheme  for affine coefficients
$\sigma(t,y,u)=\sigma_1\, u+\sigma_2$ and $b(t,x,u)=b_1\, u+b_2$ and for
$\sigma(t,y,u) = b(t,y,u)=a+b\, \, \cos(u)$. We have  studied the "experimental" speed of convergence with
respect to one mesh,  when the other one is fixed and gives a  "much
smaller" theoretical error. The second moments are computed  by Monte-Carlo
 approximations. These implementations have
been done in dimension $d=1$ for the space-time white noise $W$ and the colored noise $F$.
As expected, the observed speeds are better than the theoretical ones,  and
decrease with $\alpha$.
For example, choosing $N$ and $M$  "large" with $M\geq N^2$
  and considering "small" divisors $n$ of $N$,
we have computed  the observed linear regression coefficient  and drawn the curves of
 $\sup_{x\in [0,1]}\ln(E(|u^{n,M}(1,x)-u^{N,M}(1,x)|^2))$ as a function of
 $\ln(n)$ for various values of $\alpha$.

Note that all the results of this paper remain true if in (\ref{parablanc}) we replace
the homogeneous Dirichlet boundary conditions $u(t,x)=0$ for $x\in \partial Q$  by the homogeneous
Neumann ones $\frac{\partial u}{\partial x}(t,x)=0$ for $x\in \partial Q$.
In this last case, the eigenfunctions of $\frac{\partial}{\partial t}-\Delta$ in dimension one
is $\varphi_0(x)=1$ and for $j\geq 1$, $\varphi_j(x)=\sqrt{2}\, \cos(j\pi x)$.
Since the upper estimates of the partial sums $\sum_{j=1}^K \varphi_j(x)$ used in the
Abel transforms still hold in the case of Neumann's  conditions, the crucial
result is proved in a similar way in this case,  and the speed of convergence is preserved.
\section{Formulation of the problem}
Let $(\Omega,\FF,P)$ be a probability space, $Q=[0,1]^d$ for
some integer $d\geq 1$ and
let $F=\big(F(\varphi)\, , \, \varphi\in
\mathcal{D}(\mathbb{R}_+\times Q)\big)$ be an $L^2(P)$-valued
centered Gaussian process, which is white in time  but has a space
correlation defined as follows: given $\varphi$ and $\psi$ in
$\mathcal{D}(\mathbb{R}_+\times Q)$,
 the covariance functional of $F(\varphi)$ and $F(\psi)$ is
\begin{equation}\label{cov}
J(\varphi,\psi)=E\big( F(\varphi)\, F(\psi)\big) =
\int_0^{+\infty}\!\! dt\int\int_{(Q-Q)^*}\!\!\varphi(t,y) f(y-z) \psi(t,z) dy dz ,
\end{equation}
where  $(Q-Q)^*=\{y-z\, : y,z\in Q, \, y\neq z\}$ and
$f:(Q-Q)^*\rightarrow [0,+\infty[$ is a continuous function.
The bilinear form $J$ defined by
(\ref{cov})is non-negative definite if and only if $f$ is the
Fourier transform of a non-negative tempered distribution $\mu$ on
$Q$. Then $F$ defines a martingale-measure (still denoted by $F$),
which allows to use  stochastic integrals (see \cite{W}).
In the sequel, we suppose that for $z\in \RR^d$, $z\neq 0$,
 $f(z)=|z|^{-\alpha}$,
where $|z|$ denotes the Euclidean norm of the vector $z$. Since
$x^2+y^2\geq 2xy$, if $\alpha_j=\alpha 2^{-j}$ for $1\leq j<d$ and
$\alpha_d=\alpha 2^{-d+1}$, there exists a positive constant $C$ such that for any
 $z = (z_1, \cdots , z_d)\in\RR^d$,
\begin{equation}\label{majof} f(z)\leq C \prod_{1\leq j\leq d} f_{\alpha_j}(z_j),
\end{equation}
where $f_\alpha(\zeta)=|\zeta|^{-\alpha}$ for any $\zeta\in \RR$, $\zeta \neq 0$.
To lighten the notations, for this choice of $f$ and  $\varphi\in {\mathcal H}_d$ set
\begin{equation}\label{norm(alpha)}
\|\varphi\|_{(\alpha)}^2=\int_{{\mathbb R}^d}\int_{{\mathbb R}^d}|\varphi(y)|\, |y-z|^{-\alpha}\,
|\varphi(z)|\, dydz\, .
\end{equation}
For any $t\geq 0$, we denote by $\FF_t$ the sigma-algebra generated by
$\{ F([0,s]\times A): 0\leq s\leq t\, ,\, A\subset Q\}$.
Let $\sigma : [0,+\infty[\times Q\times \RR\rightarrow \RR$ and
$b : [0,+\infty[\times Q\times \RR\rightarrow \RR$. Suppose that there exists
a positive constant $C$ such that  for $s,t\in [0,\infty[$, $x , y \in Q$, $r,v\in \RR$,
the linear growth condition (\ref{croislin}) and 
either Lipschitz condition (\ref{lip}),
(\ref{lipn}) or (\ref{lipnm}) hold
\begin{align}\label{croislin}
&|\sigma(t,x,r)|+|b(t,x,r)|\leq C\, (1+|r|)\, ,\\
\intertext{and for $D(s,t,x,y,r,v)=|\sigma(s,x,r)-\sigma(t,y,v)|+|b(s,x,r)-b(t,y,v)|$ }
\label{lip}
&D(t,t,x,x,r,v)\leq C |r-v|,\\
\label{lipn}
& D(t,t,x,y,r,v)\leq C (|x-y|^{1-\frac{\alpha}{2}} + |r-v|),\\
\label{lipnm}
&D(s,t,x,y,r,v)\leq 
C (|t-s|^{\frac{1}{2}-\frac{\alpha}{4}}        +
|x-y|^{1-\frac{\alpha}{2}}+|r-v|).
\end{align}
For any function $u_0$ which vanishes on the boundary of $Q$,
let $u(t,x)$ denote the solution to the parabolic SPDE, which is similar to (\ref{parablanc})
\begin{equation}\label{para}
\left\{\begin{array}{l}
\frac{\partial u}{\partial t}(t,x)=\Delta u(t,x)+
\sigma(t,x,u(t,x))\frac{\partial^2 F}{\partial t\partial x}+b(t,x,u(t,x))\, ,\\
u(t,x)=0\quad \mbox{\rm for   }\; x\in\partial Q \, ,
\end{array}\right.
\end{equation}
with initial condition $u(0,x)=u_0(x)$.
Let $\NN^*$ denote the set of strictly positive integers. For any $j\in \NN^*$
and $\xi \in \RR$, set
  $\varphi_j(\xi)=\sqrt{2} \sin(j\pi \xi)$ and for
 $\underline{k}=(k_1, \cdots,\, k_d)\in \NN^{*d}$, set
\[ | \underline{k}|=\sum_{j=1}^d k_j, \quad
\varphi_{\underline{k}}(x)=\prod_{j=1}^d \varphi_{k_j}(x_{k_j})
\quad\mbox{\rm for} \quad x=(x_1,\dots, x_d)\in \RR^d .\] 
Let
$G_d(t,x,y)$ denote the Green function associated with
$\frac{\partial}{\partial t}-\Delta$ on $Q$ and homogeneous
Dirichlet boundary conditions; then  for $t>0$, $x, y\in Q$,
$G_d(t,x,y)=\sum_{\kk\in \NN^{*d} }\exp(-|\kk|^2 \pi^2 t)
\ph_{\kk}(x)\ph_{\kk}(y)$ and
\begin{equation}
\label{A1(0)} |G_d(t,x,y)|\leq C t^{-\frac{d}{2}}\,
\exp\big(-c\frac{|x-y|^2}{t}\big).
\end{equation}
When $d=1$, set $G_1=G$.  These upper estimates are classical when
the domain $Q$ has a smooth boundary under either homogeneous
Neumann's or Dirichlet's boundary conditions (see e.g. \cite{EZ},
\cite{LM}).  A simple argument shows that they can be extended for
these homogeneous conditions on the set $Q=[0,1]^d$; see e.g.
\cite{CW0} for the similar case of the parabolic operator
$\frac{\partial}{\partial t}+ \Delta^2$ on $[0,\pi]^d$. The
equation (\ref{para}) makes sense in  the following evolution
formulation (see e.g. \cite{W} for $d=1$):
\begin{align}\label{evolu}
u(t,x)=&\int_Q G_d(t,x,y)\, u_0(y)\, dy  + \int_0^t\int_Q G_d(t-s,x,y)
 \nonumber \\
& \quad \times  \left[\sigma(s,y,u(s,y))\, F(ds,dy) + b(s,y,u(s,y))\,
 ds dy\,  \right]\, .
\end{align}
We also consider the parabolic SPDE with the homogeneous boundary conditions
$\frac{\partial u}{\partial x}(t,x)=0$ for $x\in \partial Q$. Then
the functions $(\varphi_j\,  ;\,  j \geq 1)$ are replaced by $\varphi_0(\xi)=1$ and
$\varphi_j(\xi )=\sqrt{2}\, \cos(j\pi \xi)$ for $\xi \in \RR$ and $j\geq 1$.
All the other formulations remain true with $\underline{k}\in \NN^d$ instead
on $\NN^{*d}$.
\subsection{Space discretization scheme}\label{space}
As in \cite{G1}, we at first consider a finite space discretization scheme, replacing the
Laplacian by its discretization on the grid $\frac{\kk}{n}=(\frac{k_1}{n}, \cdots,
 \frac{k_d}{n})$, where $k_j\in \{0, \cdots, n\}$, $1\leq j\leq d$.
In dimension 1, we proceed as in \cite{G1}, and consider the $(n-1)\times (n-1)$-matrix
$D_n$ associated with the homogeneous Dirichlet boundary conditions and
defined by $D_n(i,i)=-2$, $D_n(i,j)=1$
if $|i-j|=1$ and $D_n(i,j)=0$ for
$|i-j|\geq 2$; then $\frac{\partial^2 u(t,x)}{\partial x^2}$ is replaced by $n^2\,
D_n \vec{u}_n(t,.)$,
where $\vec{u}^n(t)$ denotes the $(n-1)$-dimensional vector of an
 approximate solution defined on
the grid $j/n, 1\leq  j\leq n$.
In arbitrary dimension, we proceed as in \cite{CW} and
 define $D_n^{(d)}$ by induction. Let  $D_n^{(1)}=D_n$ and suppose that
 $D_n^{(d-1)}$ has been defined as a $(n-1)^{d-1}\times (n-1)^{d-1}$ matrix. Let
$Id_k$ denotes the $k\times k$ identity matrix and given a $(n-1)^{d-1}\times (n-1)^{d-1}$
matrix $A$, let $diag(A)$ denote the $(n-1)^d\times (n-1)^d$ matrix with $d-1$ diagonal
blocs equal to $A$; let $D_n^{(d)}$ denote
the $(n-1)^d\times (n-1)^d$-matrix $D_n^{(d)}$  defined by
\[D_n^{(d)}= diag(D_n^{(d-1)}) +
\left(
\begin{array}{ccccc}
-2Id_{n^{d-1}}&Id_{n^{d-1}}&0&\cdots&0\\
Id_{n^{d-1}}&-2Id_{n^{d-1}}&Id_{n^{d-1}}&\ddots&\vdots\\
0&\ddots&\ddots&\ddots&0\\
\vdots&\ddots&Id_{n^{d-1}}&-2Id_{n^{d-1}}&Id_{n^{d-1}}\\
0&\cdots&0&Id_{n^{d-1}}&-2Id_{n^{d-1}}
\end{array}\right)\, .\]
Let $\vec{u}^n(t)$ denote the $(n-1)^d$dimensional vector defined
 by $\vec{u}^n(t)_{\bf \kk}= u_n(t,{\bf \xx_{\kk}})$, with
$ {\bf \xx_{\kk}}=(x_{k_1}, \cdots, x_{k_d})$,  where
$k_j$ is the unique integer such that $ x_{k_j}=\frac{k_j-1}{n}$
 and $ k_j\in \{1,\cdots, n-1\} $
 is such that  ${\bf \kk}=(k_d-1)(n-1)^{d-1}+\cdots + (k_2-1)(n-1)+k_1.$
Let ${\mathcal L}=\big\{{\bf \xx_{\kk}}:{\bf \kk}\in \{1, \cdots, (n-1)^d\}\big\}$,
$\Box{\bf \xx_{\kk}}$ be the
lattice parallepiped of diagonal ${\bf \xx_{\kk}}=(x_{k_1}, \cdots, x_{k_d})$ and
$(x_{k_1}+\frac{1}{n}, \cdots, x_{k_d}+\frac{1}{n})$,
and set
 $F^n(t,{\bf \xx_{\kk}})= \int_{\Box{\bf \xx_{\kk}}}dF(t,x)$.
Given a function $h:[0,+\infty[\times Q\times \RR \rightarrow \RR$,
and $\vec{u}\in \RR^r$, let $h(t,x,\vec{u})=(h(t,x,u_1), \cdots,
h(t,x,u_r))$.
Then $\vec{u}^n(t)$ is solution to the following equation
\begin{equation}\label{unfaible}
d\vec{u}^n(t)=n^2\, D_n^{(d)} \vec{u}^n(t)\, dt
+ n\, \sigma(t,x,\vec{u}^n(t))\, dF(t,.), 
+b(t,x,\vec{u}^n(t)), 
\end{equation}
$\vec{u}^n(0)=\left(u_0\left(\frac{j}{n}\right),
  j\in \{1,\cdots,  n-1\}^d\right).$
We then complete $u^n(t,.)$ from the lattice ${\mathcal L}$ to $Q$ as follows.
If $d=1$, set $u^n(t,0)=u^n(t,1)=0$, $u^n(t,\frac{j}{n})=\vec{u}^n(t)_j$,
$\kappa_n(y)=[ny]/n$,
$\ph^n_j(i/n)=\ph_j(i/n)$ for $0\leq i\leq n$, and for $x\in ]i/n, (i+1)/n[$, $0\leq i<n$,
$1\leq j\leq n-1$, let
$ \ph_j^n(x)=\ph_j\left(\frac{i}{n}\right)+(nx-i)\left[ \ph_j\left(\frac{i+1}{n}\right)-
\ph_j\left(\frac{i}{n}\right)\right]$, and let 
\[ \lambda_j^n=-4\sin^2\left(\frac{j\pi}{2n}\right)n^2=-j^2\pi^2c_n^j\quad
\mbox{\rm with}\quad
c_n^j=\sin^2\left(\frac{j\pi}{2n}\right)\,
\left(\frac{j\pi}{2n}\right)^{-2} \in \left[\frac{4}{\pi^2},1\right], \] 
denote the eigenvalues of $n^2\,  D_n=n^2 \, D_n^{(1)}$; then for $t>0$, $x,y\in [0,1]$,
\begin{equation}\label{Gn}
(G_1)^n(t,x,y)=\sum_{j=1}^{n-1} \exp(\lambda_j^n\, t)\ph_j^n(x)\,
\ph_j(\kappa_n(y))\, .
\end{equation}
In dimension $d\geq 2$, we also complete the solution $u^n(t,x)$ from
$x\in {\mathcal L}$, defined as $u^n(t,{\bf \xx_k})=\vec{u}^n(t)$ to $x\in Q$
by linear interpolation, interpolating inductively on the points $(x,y)$ for
$x\in \RR^i$ and $y=(k_{i+1}/n, \cdots, k_d/n)$.
The eigenvalues and eigenvectors of $n^2\, D_n^{(d)}$ are
$ \lambda^n_{\bf \kk}=\sum_{j=1}^d \lambda^n_{k_j}$,  and
$\ph_{\bf \kk}\left(\frac{k_1\pi}{n},
\cdots, \frac{k_d\pi}{n}\right).$
For $t>0$, $x$ and
$y\in Q$ if $\kappa_n(y)=(\kappa_n(y_1), \cdots, \kappa_n(y_d))$,  let
\begin{equation}\label{Gnd}
(G_d)^n(t,x,y)=\sum_{{\bf \kk}\in \{1, \cdots,(n-1)^d\}} \exp(\lambda^n_{\bf{\kk}}t)
\ph_{\bf{\kk}}^n(x)\ph_{\bf{\kk}}^n(\kappa_n(y)).
\end{equation}
When $d=1$, simply set $G_1=G$ and $(G_1)^n=G^n$.
Then the linear interpolation of $u^n(t,.)$ from the lattice ${\mathcal L}$ to
$Q=[0,1]^d$ is solution to the evolution equation
\begin{align}\label{evolun}
&u^n(t,x)=  \int_Q (G_d)^n(t,x,y)\, u_0(\kappa_n(y))\, dy\nonumber
+ \int_0^t\int_Q (G_d)^n(t-s,x,y)\\
& \;  \times \Big[\sigma(s,\kappa_n(y), u^n(s,\kappa_n(y))F(ds,dy)
+ b(s,\kappa_n(y), u^n(s,\kappa_n(y))\, dsdy\Big]\, .
\end{align}
The $n\times n$ matrix $D_n=D_n^{(1)}$ associated with the homogeneous Neumann boundary conditions
is defined by  $D_n(1,1)=D_n(n,n)=-1$, $D_n(1,2)=D_n(n,n-1)=1$ and
for $2\leq i\leq n-1$ and $1\leq j\leq n$, $D_n(i,i)=-2$, $D_n(i,j)=1$ if $|j-i|=1$
and $D_n(i,j)=0$ for $|j-i|\geq 2$. The inductive procedure used
to construct $D_n^{(d)}$ is similar to the previous one,  replacing  1 by
$Id_{n^d}$. Then the eigenvalues of $n^2\, D_n$ are $\lambda_j^n= -4\, n^2\,
\sin^2(\frac{j\pi}{2n})=-j^2\, \pi^2\, \tilde{c}^j_n$ with $\tilde{c}^j_n\in [\frac{2}{\pi^2},1]$.
The corresponding normed eigenvectors $(e_j, 0\leq j\leq n-1)$ are again evaluations of $\varphi_j$.
More precisely, $e_j(k)=\frac{1}{\sqrt{n}}\, \varphi_j(\frac{2k-1}{2n})$ for $0\leq j\leq n-1$
and $1\leq k\leq n$. The eigenvalues $\lambda_n^{\bf \kk}$ and  the eigenfunctions $\varphi_{\bf \kk}$ of
$n^2\, D_n^{(d)}$ are defined in a way similar to the Dirichlet case, taking sums over
${\bf \kk}\in \{1,\cdots,n^d\}$. Formulas similar to (\ref{Gn}) and  (\ref{Gnd}) still hold and
(\ref{evolun}) is unchanged.
\subsection{Implicit space-time discretization scheme}\label{implicitst}
We now introduce a space-time discretization scheme. Given $T>0$, $n,m\geq 1$
we use the space mesh $1/n$ and the time mesh $T/m$,  set $t_i=iTm^{-1}$
for $0\leq i\leq m$ and  replace the time
derivative by a backward difference. Thus for $d=1$,  in the case of Dirichlet's homogeneous boundary
conditions, set $\vec{u}_0=(u_0(j/n),1\leq j\leq n-1)$
and   for $i\leq m$, set $\vec{u}_i=(u^{n,m}(iTm^{-1}, jn^{-1}), 1\leq j\leq n-1)$,
 and for $g=\sigma$ and $g=b$ let $g(t_i,.,\vec{u}_i)=
g((t_i, jn^{-1}, (u^{n,m}(t_i,jn^{-1}))\, ,\, 1\leq j\leq n-1)$. Let $\Box_{n,m}F(t_i,.)$
denote the $(n-1)$-dimensional Gaussian vector of space-time increments of $F$
on the space-time grid, i.e.,
for $1\leq j\leq n-1$, set
\begin{align*} \Box_{n,m}F(t_i,j)= nmT^{-1}& [ F(t_{i+1},(j+1)n^{-1})-F(t_i,(j+1)n^{-1})
 -F(t_{i+1},jn^{-1})+F(t_i,jn^{-1})] ;
\end{align*}
 then for every $0\leq i<m$
\begin{equation}\label{implicit1}
\vec{u}_{i+1}=\vec{u}_i+n^2\,\frac{T}{m}\, D_n\vec{u}_{i+1}+Tm^{-1}
\left[ \sigma(t_i,.,\vec{u}_i)\Box_{n,m}F(t_i,.)+b(t_i,.,\vec{u}_i)\right]\, .
\end{equation}
Since $Id-Tm^{-1}D_n$ is invertible,
\begin{align}\label{implicit2}
\vec{u}_{i+1}=&(Id-\frac{T}{m}D_n)^{-(i+1)} \vec{u}_0
+ \sum_{k=0}^i
(Id-\frac{T}{m} D_n)^{-(i-k-1)}
\left[ \sigma(t_k,.,\vec{u}_k)\Box_{n,m}F(t_k,.) +
b(t_k,.,\vec{u}_k)\right] .
\end{align}
If $d\geq 2$, we set $\Box_{n,m}F(t_i, {\bf \xx_{\kk}}  )=
n^dmT^{-1}\int_{t_i}^{t_{i+1}}\int_{\Box{\bf \xx_{\kk}}}dF(t,x)$,
and for homogeneous Dirichlet's
(resp. Neumann's) boundary conditions, define similarly
$\vec{u}_{i+1}$ as the $(n-1)^d$-dimensional (resp. $n^d$-dimen\-sio\-nal) vector such that (\ref{implicit2})
holds with $D_n^{(d)}$ instead of $D_n$. We only describe the scheme in the case
of Dirichlet's conditions;  the case of Neumann's conditions is obviously dealt with
by obvious changes.
The process $u^{n,m}$ is defined on the
space-time lattice ${\mathcal L}_T=\big\{(t_i, {\bf \xx_{\kk}}):
 0\leq i\leq m\, ,\,  {\bf \kk}\in \{1, \cdots, (n-1)^d\}\big\}$ as $(u^{n,m}(t_i, {\bf \xx_{\kk}})\, ,\,
0\leq i\leq m\, ,\,  {\bf \kk}\in \{1, \cdots, (n-1)^d\})=\vec{u}_i$~; it is then extended
to the time lattice $(t_i,x)$, $0\leq i\leq m$, $x\in Q$ as in the previous
subsection, and then extended to $[0,T]\times Q$ by time linear interpolation.
Since $\lambda_{\bf \kk}=\sum_{i=1}^d \lambda_{k_i}^n$ and $\ph_{\bf \kk}({\bf \xx_{\kk}})$
are the
eigenvalues and eigenvectors of $D_n^{(d)}$, if
\begin{equation}\label{Gnmd}
(G_d)^{n,m}(t,x,y)=\sum_{{\bf \kk}\in \{1,\cdots,(n-1)^d\}}
\left( 1-Tm^{-1}\lambda_{\bf \kk}\right)^{-[mtT^{-1}]}\; \ph_{\bf \kk}^n(x)\;
 \ph_{\bf \kk}(\kappa_n(y))\, ,
\end{equation}
then for $t=iTm^{-1}$, $1\leq i\leq m$, if for $s\in [0,T]$,
one sets $\Lambda_m(s)=[msT^{-1}]m^{-1}$ one has:
\begin{align}\label{evolimplicit}
\nonumber &u^{n,m}(t,x)=\int_Q (G_d)^{n,m}(t,x,y)u_0(\kappa_n(y)) dy +
\int_0^t\int_Q (G_d)^{n,m}(t-s+\frac{T}{m},x,y)\\
& \; \times \Big[ \sigma(\Lambda_m(s), \kappa_n(y),
 u^{n,m}(\Lambda_m(s),\kappa_n(y))) F(ds,dy) 
+ b(\Lambda_m(s), \kappa_n(y),  u^{n,m}(\Lambda_m(s),\kappa_n(y))) dy ds\Big] .
\end{align}
Again for $d=1$, let $G^{n,m}=(G_1)^{n,m}$.
\subsection{Explicit schemes}\label{explicitst}
For $T>0$, a space mesh $n^{-1}$ and a time mesh $Tm^{-1}$, we now replace
the time derivative by a forward difference. Thus if $u^n_m$ denotes the
approximating process defined for $t=t_i=iTm^{-1}$ and $x_{k_j}\in \{1, \cdots,n-1\}$,
setting  $\vec{u}_i=u^n_m(t_i,.)$, we have
\begin{equation}\label{explicit1}
\vec{u}_{i+1}=\vec{u}_i+n^2 Tm^{-1} D_n^{(d)}\vec{u}_i+Tm^{-1}
\left[ \sigma(t_i,.,\vec{u}_i)\Box_{n,m}F(t_i,.)+b(t_i,.,\vec{u}_i)\right] .
\end{equation}
In the case of homogeneous  Dirichlet boundary conditions, let $(G_d)^n_m(t,x,y)$
denote the corresponding approximation of the Green function  $G_d$ defined by
\begin{equation}\label{Gimplicit}
(G_d)^n_m(t,x,y)=\sum_{{\bf \kk}\in \{1,\cdots,(n-1)^d\}}
\left( 1+Tm^{-1}\lambda_{\bf \kk}\right)^{[mtT^{-1}]}\; \ph_{\bf \kk}^n(x)\;
 \ph_{\bf \kk}(\kappa_n(y)).
\end{equation}
Again for $d=1$, let $G_m^n=(G_1)_m^n$.
Then for $t=t_i=iTm^{-1}$, when completing the solution $u^n_m(t_i,.)$ from
the space lattice ${\mathcal L}$ to $Q$, we obtain the solution to the
following equation
\begin{align}\label{evolexplicit}
&\nonumber u^n_m(t,x)=\int_Q (G_d)^n_m(t,x,y)u_0(\kappa_n(y)) dy +
\int_0^t\int_Q (G_d)^n_m\big(t-s+T/m,x,y\big)\\
 & \quad \times \Big[ \sigma(\Lambda_m(s), \kappa_n(y),
 u^n_m(\Lambda_m(s),\kappa_n(y))) F(ds,dy) 
+ b(\Lambda_m(s), \kappa_n(y),  u^n_m(\Lambda_m(s),\kappa_n(y))) dy ds\Big] .
\end{align}
We then complete the process $u^n_m(.,x)$ by time linear interpolation and obvious changes
yield the explicit scheme for homogeneous Neumann boundary conditions.
\section{Convergence results for the discretization schemes}\label{convergence}
\setcounter{equation}{0}
In this section, we study the speed of convergence for the $d$-dimensional
space scheme and then of the $d$-dimensional implicit and explicit space-time schemes.
For the sake of simplicity, we  only write the proofs in the case of
 homogeneous Dirichlet boundary conditions.
The following result states that the solutions $u$, $u^n$, $u^{n,m}$ and $u^n_m$ exist and have
bounded moments uniformly in $n,m$. The proofs for $u$ can be found in \cite{PZ}; see also \cite{D} and  \cite{CWM}.
The arguments for the approximations are similar using (\ref{B1(1)}), (\ref{Gnalpha}),
(\ref{B1(3-1)}) and (\ref{B1(3-2)}) and the version of Gronwall's lemma in stated in \cite{G2} Lemma 3.4.
\begin{proposition}
\label{SR}
Let $u_{0}\in {\mathcal C}(Q)$ satisfy the homogeneous Neumann or Dirichlet boundary conditions,
and suppose that the coefficients $\sigma$ and $b$ satisfy the
conditions (\ref{croislin}) and (\ref{lip}); then the equation (\ref{evolu}) (resp. (\ref{evolun}), (\ref{evolimplicit})
and (\ref{evolexplicit})) has a unique solution $u$ (resp. $u^n$, $u^{n,m}$  and $u^n_m$) such that
for every $p\in [1,+\infty[$ and $T>0$:
\begin{equation}\label{solbor}
\sup_{n,m}\sup_{0\leq t\leq T}\, \sup_{x\in Q}\, \E\big(|u(t,x)|^{2p}+ |u^n(t,x)|^{2p}
+|u^{n,m}(t,x)|^{2p})+|u^n_m(t,x)|^{2p} \big) <\infty\; .
\end{equation}
\end{proposition}
We now prove H\"older regularity properties of the trajectories of $u$ and $u^n$.
 Note that for $u$,  a similar result has been proved in \cite{SS-S2} for
the heat equation with free boundary with a perturbation driven by a Gaussian process
with a more general space covariance; see also
\cite{CWM} for a related result in the case of a more general even order differential
operator.
\begin{proposition}\label{holder}
Suppose that the coefficients $b$ and $\sigma$ satisfy the Lipschitz property (\ref{lip}), that the initial condition
$u_0$ satisfies the homogeneous Dirichlet or Neumann boundary condition.\\
\indent (i) Suppose furthermore that $u_0\in  {\mathcal C}^{1-\frac{\alpha}{2}}(Q)$
and fix $T>0$.  Then, for every $p\in [1,+\infty[$,
there exists a constant $C$ such that
for $x,x'\in Q$ and $0\leq t<t'\leq T$,
\begin{align}
\sup_{0\leq t\leq T}\E\left(|u(t,x)-u(t,x')|^{2p}\right)\leq & C|x'-x|^{p(2-\alpha)}\, , \label{SR(1)}\\
\sup_{x\in Q}\E\left(|u(t',x)-u(t,x)|^{2p}\right)\leq & C|t'-t|^{p(1-\frac{\alpha}{2})}\, . \label{SR(2)}
\end{align}
\indent (ii) Suppose furthermore that $u_0\in {\mathcal C}^2(Q)$; then for every $p\in [1,+\infty[$,
there exists a constant $C$ such that for $x,x'\in Q$ and $0\leq t<t'\leq T$,
\begin{align}
\sup_{n\geq 1} \sup_{x\in Q}\E\left(|u^n(t',x)-u^n(t,x)|^{2p}\right)\leq C|t'-t|^{p(1-\frac{\alpha}{2})}\,
\label{untemps}\\
\sup_{n\geq 1} \sup_{0\leq t\leq T}\E\left(|u^n(t,x')-u^n(t,x)|^{2p}\right)\leq C|x'-x|^{p(2-\alpha)}\,
\label{unespace}
\end{align}
\end{proposition}
{\bf Proof:}  The proofs  of (\ref{SR(1)}) and (\ref{SR(2)}) can be adapted from
 Sanz-Sarr\`a \cite{SS-S2} (see also  \cite{CWM}), and are
therefore omitted. For the sake of completeness, we sketch the proof of (\ref{untemps}).
For every $t>0$, let $v^n(t,x)=\int_Q (G_d)^n(t,x,y)\, u_0(\kappa_n(y)\, dy$
and $w^n(t,x)=u^n(t,x)-v^n(t,x)$, where  $(G_d)^n$ is the fundamental solution
of   $\frac{\partial }{\partial t} -\Delta_n=0$,
\begin{align}\label{Deltan}
\Delta_n U(y)=&n^2 \sum_{i=1}^d \Big[U\Big(\sum_{j\neq i}\frac{[ny_j]}{n} e_j + \frac{[ny_i]+1}{n} e_i \Big)
-2 U\Big(\frac{[ny]}{n}\Big)
+ U\Big(\sum_{j\neq i}\frac{[ny_j]}{n} e_j+ \frac{[ny_i]-1}{n} e_i\Big)\Big]  ,
\end{align}
and $(e_i, 1\leq i\leq d)$ denotes the canonical basis of $\RR^d$. 
Then if $u_0\in {\mathcal C}^2(Q)$,
\begin{equation}\label{evolvn}
v^n(t,x)=u_0(x)+\int_0^t\int_Q (G_d)^n(s,x,y)\, \Delta_n u_0(y)\, dy\,ds\, .
\end{equation}
Using the fact that $\Delta_n u_0$ is bounded if $u_0\in {\mathcal C}^2(Q)$, and (\ref{*1}),
we deduce that for any $\lambda >0$
\begin{equation}\label{vn-vn}
\sup_{n\geq 1}\sup_{x\in Q} |v_n(t,x)-v_n(t',x)| 
\leq C\, |t'-t|^{1-\lambda}\, .
\end{equation}
Computations similar to those used in \cite{SS-S2}, using Burkholder's and H\"older's inequalities with
respect to suitable measures,
 (\ref{*1})-(\ref{*4}) and (\ref{solbor}), show the existence
of $C_p>0$ such that for any $0\leq t<t'\leq T$
\begin{equation}\label{wn-wn}
\sup_{n\geq 1}\sup_{x\in Q} \E(|w^n(t,x)-w^n(t',x)|^{2p})\leq C_p\, |t'-t|^{p(1-\frac{\alpha}{2})}\, .
\end{equation}
The inequalities (\ref{vn-vn}) and (\ref{wn-wn}) conclude the proof of (\ref{untemps}).
\qquad $\Box$

The first convergence  result of this section is that of $u^n$ to $u$.
\begin{Th}\label{appron}
Let $\sigma$ and $b$ satisfy the conditions (\ref{croislin}) and (\ref{lipn}), $u$ and $u^n$ be the solutions
to (\ref{evolu}) and (\ref{evolun}) respectively, where the Green functions $G_d$ and $(G_d)^n$ are defined with
the homogeneous  Neumann or Dirichlet boundary conditions on $Q$.\\
(i) If the initial condition $u_0$ belongs to $\mathcal{C}^3(Q)$,
 then for every $T>0$ and $p\in [1+\infty[$,  there exists a constant $C_p(T)>0$ such that:
\begin{equation}\label{vitesnc3}
\sup_{(t,x)\in[0,T]\times Q}\; \E\left( \, |u(t,x)-u^n(t,x)|^{2p}|\, \right)\leq C_p(T)\, n^{-(2-\alpha)\, p}\, .
\end{equation}
(ii) If the initial condition $u_0$ belongs to  $\mathcal{C}^{1-\frac{\alpha}{2}}(Q)$, then there exists $\nu >0$ such
that given any  $p\in [1,+\infty[$, there exists a constant $C_p>0$ such that,  for every $t>0$:
\begin{equation}\label{vitesnc1-alpha}
 \sup_{x\in  Q}\; \E\left( \, |u(t,x)-u^n(t,x)|^{2p}|\, \right)\leq C_p\,t^{-\nu}\,  n^{-(2-\alpha)\, p}\, .
\end{equation}
(iii) Finally, if $u_0$ belongs to $\mathcal{C}_0(Q)$, then for
all $p\in [1,+\infty[$, as $n\rightarrow +\infty$, 
${\displaystyle \sup_{(t,x)\in [0,T]\times Q}\! \E\big( |u(t,x)}$ \linebreak
 ${- u^n(t,x)|^{2p}|
\big)}$ converges to 0, and the sequence $u^n(t,x)$ converges a.s.
to $u(t,x)$ uniformly on $[0,T]\times Q$.
\end{Th}
{\bf Proof:} 
  As in \cite{G1}, set $v(t,x)=\int_Q G_d(t,x,y)\, u_0(y)\, dy$,
 $v^n(t,x)=\int_Q (G_d)^n(t,x,y)\, u_0(\kappa_n(y))\, dy$, and $u(t,x)=v(t,x)+w(t,x)$
 $u^n(t,x)=v^n(t,x)+w^n(t,x)$.
If $u_0\in \mathcal{C}^{1-\frac{\alpha}{2}}(Q)$ (and hence is bounded), using 
(\ref{B2(0)}) and (\ref{B1(1)}), we deduce that
for any $\lambda\in ]0,1[$, there exists $\mu >0$, $C>0$ such that for $t>0$, $\nu=\lambda\vee \mu$,
\begin{align}\label{v-vn0}
& \sup_{x\in Q} |v(t,x)-v^n(t,x)|\, \leq\,  \int_Q\, \Big[\, |G_d(t,x,y)-(G_d)^n(t,x,y)|\, |u_0(y)| \nonumber \\
&\qquad + |(G_d)^nt,x,y)|\,
|u_0(y)-u_0(\kappa_n(y))|\,\Big]\,  dy  
\leq\,   
C\, (1+t^{-\nu})\,e^{-ct}\,  n^{-(1-\frac{\alpha}{2})}\, .
\end{align}
If $u_0\in \mathcal{C}^3(Q)$, then since $G_d$ (resp. $(G_d)^n$) is the fundamental solution
of $\frac{\partial }{\partial t} -\Delta=0$ (resp  $\frac{\partial }{\partial t} -\Delta_n=0$),
where  $\Delta_n$ is defined by (\ref{Deltan}), integrating by parts we deduce that
$v(t,x)=u_0(x)+\int_0^t\int_Q G_d(s,x,y)\, \Delta u_0(y)\, dy$ and
$v^n(t,x)=u_0(x)+\int_0^t\int_Q (G_d)^n(s,x,y)\, \Delta_n u_0(y)\, dy$.
Hence $|v(t,x)-v^n(t,x)|\leq \sum_{i=1}^3 A_i(t,x)$, where
 \begin{align*}
A_1(t,x)&=|u_0(t,x)-u_0(\kappa_n(x))|\, ,    \quad 
A_2(t,x)=\Big|\int_0^t\int_Q [G_d(s,x,y)-(G_d)^n(s,x,y)] \Delta u_0(y) dy ds\Big|,\\
A_3(t,x)&=\Big|\int_0^t\int_Q (G_d)^n(s,x,y) [\Delta u_0(y)- \Delta_n u_0(\kappa_n(y))] dy ds\Big|.
\end{align*}
Since $\Delta u_0$ is bounded and $\|u_0(.)-u_0(\kappa_n(.))\|_\infty +
\| \Delta u_0(.) -\Delta_n u_0(\kappa_n(.))\|_\infty\leq C\, n^{-1}$, the inequalities
(\ref{B1(1)}) and (\ref{B2(1)}) imply
\begin{equation}\label{v-vn1}
\sup_{(t,x)\in [0,+\infty[\times Q} |v(t,x)-v^n(t,x)|\leq C\, n^{-1}\, .
\end{equation}
Furthermore, for every  $t\in ]0, T]$, $ \sup_{x\in Q} \E(|w(t,x)-w^n(t,x)|^{2p})\leq C
 \sum_{i=1}^6 B_i(t)$,
where
\begin{align}\label{defB}
B_1(t)=& \sup_{x\in Q} \E\Big( \Big| \int_0^t\!\! \int_Q\!\! G_d(t-s,x,y) \Big(\sigma(s,y,u(s,y))
-\sigma(s,\kappa_n(y), u(s,\kappa_n(y))\Big)
 F(ds,dy)\Big|^{2p} \Big)  ,\nonumber\\
B_2(t)=& \sup_{x\in Q} \E\Big( \Big| \int_0^t\!\!\! \int_Q\!\!\!  G_d(t-s,x,y)  \Big(
\sigma(s,\kappa_n(y),u(s,\kappa_n(y)))
\nonumber \\
&\qquad\qquad \qquad
- \sigma(s,\kappa_n(y), u^n(s,\kappa_n(y))\Big)   F(ds,dy)\Big|^{2p} \Big) ,\nonumber \\
B_3(t)=& \sup_{x\in Q} \E\Big( \Big| \int_0^t\! \! \int_Q \!\! G_d(t-s,x,y)-(G_d)^n(t-s,x,y)\Big)
  \sigma(s,\kappa_n(y), u^n(s,\kappa_n(y))
 F(ds,dy)\Big|^{2p} \Big) , \nonumber \\
B_4(t)=& \sup_{x\in Q} \E\Big( \Big| \int_0^t\!\! \int_Q G_d(t-s,x,y)  \Big( b(s,y,u(s,y))
 -b(s,\kappa_n(y), u(s,\kappa_n(y))\Big)
 dy ds\Big|^{2p}\, \Big) ,\nonumber \\
 B_5(t)=& \sup_{x\in Q} \E\Big( \Big| \int_0^t\!\! \int_Q \!\! G_d(t-s,x,y)   
 \Big( b(s,\kappa_n(y),u(s,\kappa_n(y)))
-  b(s,\kappa_n(y), u^n(s,\kappa_n(y))\Big)   dy ds\Big|^{2p} \Big) ,\nonumber \\
 B_6(t)=& \sup_{x\in Q} \E\Big(\, \Big| \int_0^t\!\! \int_Q \!\!\Big( G_d(t-s,x,y)-(G_d)^n(t-s,x,y)\Big)
b(s,\kappa_n(y), u^n(s,\kappa_n(y))
 dy ds \Big|^{2p} \Big).
\end{align}
Burkholder's inequality, (\ref{A1(1)}), H\"older's inequality with respect to
$\| G_d(t-s,x,.)\|_{(\alpha)}^2\, ds$, Fubini's theorem, (\ref{lipn}),
Schwarz's inequalities and (\ref{SR(1)}) imply that
\begin{align}\label{majoB1n}
&B_1(t)\leq  C_p \sup_{x\in Q} \E  \Big|\! \int_0^t
\left\| |G_d(t-s,x,.)| \big( n^{-\frac{2-\alpha}{2}}+ |u(s,.)-u(s,\kappa_n(.)| \big)\right\|^2_{(\alpha )}
 ds \Big|^p \nonumber \\
 & \qquad \leq  C_p\, 
\!  \int_0^t (t-s)^{-\frac{\alpha}{2}}
\! \Big[\,  n^{-p\, (2-\alpha)} +\! \!\sup_{(s,\xi)\in [0,t]\times Q}\!\!\!
\E(|u(s,\xi)-u(s,\kappa_n(\xi))|^{2p}) \Big]
 ds 
 \leq C_p\, n^{-p(2-\alpha)} .
\end{align}
Similar arguments based on (\ref{A1(1)}) and  (\ref{lipn}) (resp. 
(\ref{B2(2)}), (\ref{croislin}) and (\ref{solbor})) imply  that 
\begin{align}\label{majoB2n}
B_2(t)\leq& C_p 
\int_0^t (t-s)^{-\frac{\alpha}{2}} \Big( \sup_{x\in Q} |v(s,x)-v^n(s,x)|^{2p} 
+ \sup_{x\in Q} \E(|w(s,x)-w^n(s,x)|^{2p})\Big) ds  ,\\
\label{majoB3n}
  B_3(t) \leq & C_p\, n^{-(2-\alpha)\, p}\, .
\end{align}
The deterministic integrals are easier to deal with. Using H\"older's inequality with respect to the
measure $|G(t-s,x,y)|\, dy\, ds$, (\ref{A1(0)}), (\ref{lipn}) and (\ref{SR(1)}) we deduce that
\begin{align}\label{majoB4n}
 B_4(t)\leq & C_p \sup_{x\in Q} \E \Big( \Big| \int_0^t\int_Q
|G_d(t-s,x,y)| dy ds\Big)^{2p-1}
\int_0^t\!\int_Q |G_d(t-s,x,y)|\nonumber \\
&   \times  \Big(n^{-(2-\alpha)\, p} +
 \E(|u(s,y)-u(s,\kappa_n(y))|^{2p}) \Big) dy dz ds
\leq  C n^{-p(2-\alpha)} .
\end{align}
Similarly, H\"older's  inequality,  (\ref{lipn}) and (\ref{A1(0)}) (resp. 
(\ref{croislin}), (\ref{solbor}) and  (\ref{B2(1)})) imply
\begin{align}\label{majoB5n}
B_5(t)
\leq & C \int_0^t\! \Big[ \sup_{x\in Q} |v(s,x)-v^n(s,x)|^{2p} +
 \sup_{x\in Q} \E(|w(s,x)-w^n(s,x)|^{2p}\Big] ds ,\\
\label{majoB6n}
 B_6(t) 
\leq & C\, n^{-2p}\, .
\end{align}
The inequalities (\ref{majoB1n})-(\ref{majoB6n}) imply that for any $T>0$ and $p\in [1,+\infty[$,
there exists a constant
$C>0$ such that for $0\leq t\leq T$,
\begin{align}\label{w-wn}
&\sup_{x\in Q} \E(|w(t,x)-w^n(t,x)|^{2p})\leq C\Big(  n^{-p(2-\alpha)} + \int_0^t (t-s)^{-\frac{\alpha}{2}} \nonumber \\
&\qquad\times \Big[
\sup_{x\in Q}  |v(s,x)-v^n(s,x)|^{2p} + \sup_{x\in Q} 
\E(|w(s,x)-w^n(s,x)|^{2p}\big)\Big]  ds \Big)  .
\end{align}
Thus, (\ref{v-vn1}) and  Gronwall's lemma (see e.g. \cite{G2}, lemma 3.4) imply that if $u_0\in\mathcal{C}^3(Q)$,
\begin{align*}
 \sup_{x\in Q}&  \E(|w(t,x)-w^n(t,x)|^{2p})\leq C_p\, \Big[ n^{-p (2-\alpha)} + \int_0^t (t-s)^{-\frac{\alpha}{2}}\\
& \qquad \times \sup_{x\in Q} \E(|w(s,x)-w^n(s,x)|^{2p}) \, ds \Big]
\leq C_p\, n^{-p(2-\alpha)}.
\end{align*}
This inequality together with (\ref{v-vn1}) yield (\ref{vitesnc3}). If $u\in \mathcal{C}^{1-\frac{\alpha}{2}}(Q)$,
 using again Gronwall's lemma and (\ref{v-vn0}), we deduce that for some $\lambda \in ]0,1[$, one has
\begin{align*}
&\sup_{x\in Q} \E(|w(t,x)-w^n(t,x)|^{2p})\leq C_p \Big( n^{-p (2-\alpha)}
+ \int_0^t (t-s)^{-\frac{\alpha}{2}} \Big[  s^{-\lambda}  n^{-p (2-\alpha)}  \\
& \qquad\qquad  \qquad + 
\sup_{x\in Q} \E(|w(s,x)-w^n(s,x)|^{2p})\Big] ds \Big)\leq C_p n^{-p(2-\alpha)}.
\end{align*}
This inequality and (\ref{v-vn0}) imply (\ref{vitesnc1-alpha}).

Finally, let $u_0\in \mathcal{C}^0(Q)$ and  for any $\varepsilon >0$,
 let $u_{0,\varepsilon}$ denote a function in $\mathcal{C}^3(Q)$
such that $\|u_0-u_{0,\varepsilon}\|_\infty \leq \varepsilon$. Let $u_\varepsilon = v_\varepsilon + w_\varepsilon$ and
$u_\varepsilon^n=v_\varepsilon^n+w_\varepsilon^n$ denote the previous decompositions of the solution $u_\varepsilon$ and
its space discretization $u_\varepsilon^n$ with the initial condition $u_{0,\varepsilon}$.
Then
\begin{align}\label{v-vn2}
&\sup_{(t,x)\in [0,T]\times Q} |v(t,x)-v^n(t,x)| \leq
\sup_{(t,x)\in [0,T]\times Q}|v_\varepsilon(t,x)-v_\varepsilon^n(t,x)|
\nonumber \\
&\qquad +\Big| \int_Q G_d(t,x,y) |u_0(y)-u_{0,\varepsilon}(y)] dy\Big|
  + \Big| \int_Q (G_d)^n(t,x,y) [u_0(\kappa_n(y))-u_{0,\varepsilon}(\kappa_n(y))] dy\Big|
\nonumber \\
&\quad \leq C\varepsilon + \sup_{(t,x)\in [0,T]\times Q}|v_\varepsilon(t,x)-v_\varepsilon^n(t,x)|.
\end{align}
Hence (\ref{w-wn}) and (\ref{v-vn2}) imply that
\[ \sup_{x\in Q}   \E(|w(t,x)-w^n(t,x)|^{2p}) \leq C
 \Big[ \varepsilon + n^{-p(2-\alpha)}
 + \int_0^t (t-s)^{-\frac{\alpha}{2}}
\sup_{x\in Q} \E(|w(s,x)-w^n(t,x)|^{2p}) ds \Big].
\] 
Gronwall's lemma concludes the proof of the theorem.
\qquad $\Box$

 We now prove the convergence of $u^{n,m}$  and of $u^n_m$ to $u^n$ as $m\rightarrow +\infty$.
\begin{Th}\label{n-nm}
Let $\sigma$ and $b$ satisfy the conditions (\ref{croislin}) and (\ref{lipnm}). Then\\
(i) If $u_0\in {\mathcal C}^2(Q)$, then for every $T>0$ and $p\in [1,+\infty[$, there exists a constant $C_p(T)>0$ such
that
\begin{equation}\label{convn-nm}
\sup_{n\geq 1}\, \sup_{t\in [0,T]}\, \sup_{x\in Q} \E(|u^n(t,x)-u^{n,m}(t,x)|^{2p}) \leq C_p(T)\, m^{-p (1-\frac{\alpha}{2})}\, .
\end{equation}
(ii) If $u_0\in {\mathcal C}(Q)$, then  $\sup_{n\geq 1}\, \sup_{t\in [0,T]}\,
\sup_{x\in Q} |u^n(t,x)-u^{n,m}(t,x)|$ converges to 0 as $m\rightarrow +\infty$ and for every
$t>0$ and $p\in [1,+\infty[$ there exists a constant $C_p(t)$ such that
\[ \sup_{n\geq 1}\, \sup_{x\in Q} \E(|u^n(t,x)-u^{n,m}(t,x)|^{2p}) \leq C_p(t)\, m^{-p (1-\frac{\alpha}{2})}\, .\]
(iii) The results of (i) and (ii) hold with $u^n_m$ instead of $u^{n,m}$ if one requires that
 $\frac{n^2\, T}{m}\leq q<\frac{1}{2}$.
\end{Th}
{\bf Proof:} 
Let $v^n(t,x)=\int_Q (G_d)^n(t,x,y) u_0(\kappa_n(y))dy$,
$v^{n,m}(t,x)=\int_Q (G_d)^{n,m}(t,x,y) u_0(\kappa_n(y)) dy$.
Suppose at first that $u_0\in {\mathcal C}^2(Q)$ and as in the
proof of (3.23) in \cite{G2},  for $d=1$ set   $ I=\sup_{t\in
[0,T]}\sup_{x\in [0,1]}  |v^{n,m}(t,x)
  -v^n(t,x)|
\leq \sum_{i=1}^3
I_i$, where
\begin{align*}
I_1=&\sup_{t\in [0,T]}\sup_{x\in [0,1]} \big|v^{n,m}\big( \big[{mt}{T}^{-1}\big]\, {T}{m}^{-1},x\big)
- v^n\big( \big[{mt}{T}^{-1}\big]\, {T}{m}^{-1},x\big)\big| , \\
I_2=&\sup_{t\in [0,T]}\sup_{x\in [0,1]} \big|v^n\big( \big[{mt}{T}^{-1}\big]\, {T}{m}^{-1},x\big)
- v^n(t,x)\big|, \\
I_3=&\sup_{t\in [0,T]}\sup_{x\in [0,1]} \big|v^n \big( \big[\big({mt}{T}^{-1}+1\big)\big]\, {T}{m}^{-1},x\big)
-v^n(t,x)\big|\, .
\end{align*}
The inequalities (3.27) and (3.28) in \cite{G2} imply that $I_2+I_3 \leq C\, m^{-\frac{1}{2}}$. Furthermore, using
an estimate of \cite{G2}, we deduce that
\begin{align*}
I_1\leq &\,  C\, \sup_{n\geq 1}\, \sup_{t\in [0,T]}\, \sup_{x\in [0,1]} \sum_{j=1}^{n-1} j^{-2}\,
\exp\Big( \lambda_j^n\, \big[\frac{mt}{T}\big]\, \frac{T}{m}\Big)
\Big| 1 - \exp\Big[ \big[\frac{mt}{T}\big]\,
\Big( \lambda_j^n \, \frac{T}{m}+\ln\Big( 1-\lambda_j^n\,  \frac{T}{m}\Big)\Big] \Big|\, .
\end{align*}
For $t\leq T\, m^{-1}$, $[\frac{mt}{T}]=0$ and the right hand-side of the previous inequality is 0.
If $t\geq T\, m^{-1}$, then there exists a constant $c>0$ such that $\frac{T}{m}\, [\frac{mt}{T}] \geq c\, t$ and using
(\ref{A4(1)}) we deduce that
\begin{align*}
I_1\leq&C\, \sup_{n\geq 1}\sup_{t\in [\frac{T}{m},T]}\,  \sum_{j=1}^{n-1} j^{-2}\, e^{-ctj^2}\,
\big| 1-\exp(-j^4\, t\, m^{-1})\big|
\\
\leq & \sup_{n\geq 1}\sup_{t\in [\frac{T}{m},T]}\, m^{-1}\,   \sum_{j=1}^{n-1} j^2\, t \, e^{-ctj^2}\,
\leq C\, \sup_{n\geq 1}\sup_{t\in [\frac{T}{m},T]}\, m^{-1}\,   \sum_{j=1}^{n-1} e^{-ctj^2} \leq C\, m^{-\frac{1}{2}}\, .
\end{align*}
Hence for $d=1$,
$  \sup_{n\geq 1}\, \sup_{t\in [0,T]}\, \sup_{x\in Q}\, |v^n(t,x)-v^{n,m}(t,x)|\leq C\, m^{-\frac{1}{2}}$,
and an easy argument shows that this inequality can be extended to any $d\geq 1$. Furthermore, for any $m \geq 1$
and $t\in [0,T]$,
$\sup_{n\geq 1} \sup_{x\in Q}\E(|w^n(t,x) - w^{n,m}(t,x)|^{2p})\leq C\, \sum_{i=1}^6 \tilde{B}_i(t)$, where
 $\tilde{B}_1(t)$ and $\tilde{B}_4(t)$ are similar to  ${B}_1(t)$ and ${B}_4(t)$
in the proof of (\ref{defB}) replacing $\varphi(s,y,u(s,y))-
 \varphi(s,\kappa_n(y),u(s,\kappa_n(y))) $ by  $\varphi(s,\kappa_n(y),u^n(s,\kappa_n(y)))
 -\varphi(\Lambda_m(s),\kappa_n(y), u^n(s,\kappa_n(y)))$, 
$\tilde{B}_2(t)$ and $\tilde{B}_5(t)$ are similar to  ${B}_2(t)$ and ${B}_5(t)$ replacing
$ \varphi(s,\kappa_n(y),u(s,\kappa_n(y)))
 - \varphi(s, 
\kappa_n(y), u^{n}(s,\kappa_n(y)))$
by $\varphi(\Lambda_m(s),\kappa_n(y),u^n(\Lambda_m(s),\kappa_n(y)))
 - \varphi(\Lambda_m(s),\kappa_n(y), 
u^{n,m}(\Lambda_m(s),\kappa_n(y)))$
with $\varphi = \sigma$ or $b$ respectively, and finally $\tilde{B}_3(t)$ and $\tilde{B}_6(t)$
are similar to  ${B}_3(t)$ and ${B}_6(t)$ replacing $G_d-(G_d)^n$ by $(G_d)^n-(G_d)^{n,m}$.
The argument is then  similar to that used in the proof of Theorem \ref{appron}.
The inequalities (\ref{lipnm}), (\ref{B1(3-1)}), (\ref{solbor})
and (\ref{untemps}) provide   an upper estimate of $\tilde{B}_1$,  (\ref{B3(4)}) and (\ref{solbor})
 give an upper estimate of $\tilde{B}_3$
so that  $\tilde{B}_1(t)+\tilde{B}_3(t)\leq C\, m^{-(1-\frac{\alpha}{2})\, p}$.
On the other hand, (\ref{B1(3-2)}) and (\ref{lipnm})
show that for some $\lambda\in ]0,1[$,
\[\tilde{B}_2(t)\leq \int_0^t (t-s)^{-\lambda}\, \sup_{n\geq 1}\, \sup_{y\in Q}\, \E\big(|u^n(\Lambda_m(s),\kappa_n(y))
- u^{n,m}(\Lambda_m(s),\kappa_n(y))|^{2p}\,  \big)\, ds \, .\]
A similar argument based on (\ref{B1(3-1)}), (\ref{B3(3)}), (\ref{solbor})  (\ref{untemps}) 
proves that 
$\tilde{B}_4(t)+\tilde{B}_6(t)\leq C\, m^{-\mu}$ for any $\mu\in ]0,1[$ and shows that for some $\lambda\in ]0,1[$,
\[\tilde{B}_5(t)\leq \int_0^t (t-s)^{-\lambda}\, \sup_{n\geq 1}\, \sup_{y\in Q}\, \E\big(|u^n(\Lambda_m(s),\kappa_n(y))
- u^{n,m}(\Lambda_m(s),\kappa_n(y))|^{2p}\, \big) \, ds\, .\]
Thus, Gronwall's lemma concludes the proof of (\ref{convn-nm}).
The rest of the proof of the theorem, which is similar to that of Theorem \ref{appron} is omitted. \qquad $\Box$
\section{Refined estimates of differences of Green functions } \label{cruciale}
This section is devoted to prove some crucial evaluations for the norms of the difference  between  $G_d$ and its
space discretizations $(G_d)^n$ , $(G_d)^{n,m}$ or $(G_d)^n_m$; indeed, as shown in the previous section, they
provide the speed of convergence of the scheme. We suppose again that these kernels are defined
in terms of the homogeneous Dirichlet boundary conditions. Simple modifications of the proof yield  similar 
estimates for the homogeneous Neumann ones.

The main ingredient in the proofs will be the so-called {\it Abel's summation method}, which is a discrete "integration-by-parts"
formula and is classically used in analysis to
evaluate non absolutely convergent series. More precisely :
\begin{quote}
\it Let $(a_{n})_{n\in\NN}$, $(b_{n})_{n\in\NN}$ be sequences of real numbers, $A_{-1}=0$
and   $A_{n}=\sum_{k=0}^na_{k}$ if $n\geq 0$. 
Then, for any $0\leq N_{0}< N$, one has  
\[ \sum_{k=N_{0}}^N a_{k}\, b_{k}=\sum_{k=N_{0}}^N(A_{k}-A_{k-1})b_{k}
=A_{N}b_{N}-A_{N_{0}-1}b_{N_{0}}-\sum_{k=N_{0}}^{N-1}A_{k}\, (b_{k}-b_{k+1})\, . \]
\end{quote}
In particular, this technique will be employed repeatedly throughout the proofs with
$x\in ]0,2[$ and  $a_{k}=\cos(k\pi x)$,
for which the corresponding sequence $A_{k}$ satisfies the property
$|A_{k}| \leq\frac{C}{\left|\sin\left(\frac{\pi x}{2}\right)\right|}$,
or $a_{k}=\sin(k\pi x)$, for which $A_k$ satisfies a similar inequality,
and various monotonous
sequences $( b_{k})$;  see \cite{abel} pages 17 - 18
for a more detailed account on the subject.
\begin{lemme}
\label{B2}
There exists  positive  constants $c$,
  $C$,  $\mu$  
such that for $t>0$,  
  $n\geq 2$:
\begin{align}
&\sup_{x\in Q} \|G_d(t,x,.)-(G_d)^n(t,x,.)\|_1\leq  C\, n^{-1}\, (1+t^{-\mu})\, e^{-ct}\, ,
\label{B2(0)}\\
&\int_{0}^{+\infty}\hspace{-2mm}\sup_{x\in Q} \int_Q |G_d(t,x,y)-(G_d)^n(t,x,y)|\,dydt\leq Cn^{-1}\, ,
\label{B2(1)}\\
&\int_{0}^{+\infty}\sup_{x\in Q}\|\, G_d(t,x,\cdot)-(G_d)^n(t,x,\cdot) \|_{(\alpha)}^2\, dt
\leq C n^{-(2-\alpha)}\, .
\label{B2(2)}
\end{align}
\end{lemme}
{\bf Proof : }Let $\gamma >0$ to be fixed later on; the inequalities 
(\ref{B1(1)}), (\ref{A1(1)})
and
(\ref{B1(2)}) imply that for $0<\lambda<1$,
\begin{align}
& \int_0^{\gamma n^{-2}}\!\!\sup_{x\in Q} \|G_d(t,x,.)-(G_d)^n(t,x,.)\|_1\, dt
\leq C\, n^{-2+\lambda}
\, ,\label{B2(3)}\\
& \int_0^{\gamma n^{-2}}\!\!\sup_{x\in Q} \|\, G_d(t,x,.)-(G_d)^n(t,x,.) \|_{(\alpha)}^2 \, dt
\leq C\, n^{-2+\alpha}\, .\label{B2(4)}
\end{align}
To 
estimate $\int_{\gamma\, n^{-2}}^{+\infty}\sup_x \| G_d(t,x,.) -
(G_d)^n(t,x,.) \| \, dt$, where $\|\; \|$ denotes either the $\|
\; \|_1$ or $\|\;  \|_{(\alpha)}$ norm, we  first deal with the
case $d=1$ and $\alpha <1$.

\noindent{\bf Case $d=1$ and $\alpha <1$.}  As in Gy\"ongy \cite{G1},  write
$|G(t,x,y)-G^n(t,x,y)|\leq \sum_{i=1}^4T_{i}(t,x,y)$,
where
\begin{align}
T_{1}(t,x,y)=&\Big|\sum_{j\geq n} e^{-j^2\pi^2t}\; \varphi_{j}(x)\; \varphi_{j}(y)\Big|\, ,\nonumber\\
T_{2}(t,x,y)=&\Big|\sum_{1\leq j \leq n-1}\left[e^{\lambda_{j}^nt}-e^{-j^2\pi^2t}\right]\;
\varphi_{j}(x)\; \varphi_{j}(y)\Big|\, ,\nonumber\\
T_{3}(t,x,y)=&\Big|\sum_{1\leq j\leq n-1}e^{\lambda_{j}^nt}\;
 [\varphi_{j}(x)-\varphi_{j}^n(x)]\; \varphi_{j}(y)\Big|\, ,\nonumber\\
T_{4}(t,x,y)=&\Big|\sum_{1\leq j\leq n-1}e^{\lambda_{j}^nt}\; \varphi_{j}^n(x)\;
 [\varphi_{j}(y)-\varphi_{j}(k_{n}(y))]\Big|\, . 
\label{B2(5)}
\end{align}
To study the $\| \|_{(\alpha)}$ norm of a non-negative function $R(t,x,.)$, for $x\in [0,1]$
 and $i\in \{1,2,3\}$, let
\begin{equation}\label{B2(7)}
A_n^i(x)=\{y\in [0,1]: |y-x|\leq i\, n^{-1} \; \mbox{\rm or }\; y+x\leq i\,
 n^{-1} \; \mbox{\rm or }\;
2-x-y\leq i\, n^{-1}\,  \}\, .
\end{equation}
Then $dy(A_n^i(x))\leq Cn^{-1}$ and for $x\in [0,1]$, $y,z\in A_n^i(x)$,
 $|y-z|\leq 2i\, n^{-1}$; furthermore, 
\[ \|R(t,x,.)\|_{(\alpha)}^2\leq 2\, \left[\,\|R(t,x,.)\, 1_{A_n^2(x)}(.)\|^2_{(\alpha)} +
\|R(t,x,.)\, 1_{A_n^2(x)^c}(.)\|^2_{(\alpha)}\, \right]\, .\]
Set
${\mathcal A}^{(1)}_n(x)=\{ (y,z)\in Q^2 : |y-x|\vee|z-x|\leq 2\, n^{-1}\}$,  
${\mathcal A}^{(2)}_n(x)=\{ (y,z)\in Q^2 : |y-x|\vee (x+z)\leq 2\, n^{-1}\},$ and  
${\mathcal A}^{(3)}_n(x)=\{ (y,z)\in Q^2 : |y-x|\vee (2-x-z)\leq 2\, n^{-1}\}$
and for $i=1,2,3$, let
\begin{equation}\label{mathcalA}
R^{(i)}(t,x)=\int_{{\mathcal A}^{(1)}_n(x)}\!\!\! R(t,x,y) |y-z|^{-\alpha}
 R(t,x,z)\, dydz  \,. 
\end{equation}
Then $\|R(t,x,.) 1_{A_n^2(x)}(.)\|^2_{(\alpha)} \leq C  \sum_{i=1}^3  R^{(i)}(t,x) $.
Let ${\mathcal B}^{(1)}_n(x)=\{ (y,z)\in Q^2\, :\, 2\, n^{-1}\leq |y-x|\wedge |z-x| \, ,\, |y-z|\leq 2\, n^{-1}\}$ and
${\mathcal B}^{(2)}_n(x)=\{ (y,z)\in Q^2\, :\, 2\, n^{-1}\leq |y-z|\wedge
|y-x|\wedge |z-x|\}$ and  for $i=1,2$ set
\begin{equation}\label{mathcalB}
\bar{R}^{(i)}(t,x)=\int_{{\mathcal B}^{(i)}_n(x)}\!\!\!  R(t,x,y) |y-z|^{-\alpha}
R(t,x,z)\, dy dz\, .
\end{equation}
Then $ \| R(t,x,.) 1_{A_n^2(x)^c}(.) \|_{(\alpha )}^2
 \leq C \sum_{i=1}^2 \bar{R}^{(i)}(t,x)$.
These notations will be used repeatedly throughout the proof for various functions $R$.

\noindent {\bf Estimate of  $T_{2}$}
 This  term is the most delicate  to handle.
 Set $\Delta_{j}^n(t):=e^{\lambda_{j}^nt}-e^{-j^2\pi^2t}$. Then
 for any $A\in [0,2]$ we have
\begin{equation}
\label{B2(10)}
0\leq \Delta_{j}^n(t)\leq C\,  (j/n)^2\, j^2\, t\, e^{-cj^2t}
\leq C\, n^{-A}\, j^A\, e^{-cj^2t}\, ,
\end{equation}
so that (\ref{A4(1)}) with  $K=A$ yields
\begin{equation}\label{B2(11)}
\sup_{x\in [0,1]} T_2(t,x,y) \leq C\, n^{-A}\, t^{-\frac{A+1}{2}}\, e^{-ct}\, .
\end{equation}
Furthermore,  $T_2(t,x,y)\leq |\sum_{j=1}^{n-1} e^{\lambda^n_jt}\, \varphi_j(x)\, \varphi_j(y)|+
|\sum_{j=1}^{n-1} e^{-j^2\pi^2 t}\, \varphi_j(x)\, \varphi_j(y)|$, and 
  Abel's summation method
yields that for  $1\leq N_1(n)<N_2(n)\leq n-1$,
\begin{equation}\label{B2(12)}
\Big|\sum_{j=N_1(n)}^{N_2(n)} \Delta_j^n(t)\, \varphi_j(x)\, \varphi_j(y)\Big| \leq  e^{-c\, N_1(n)^2\, t}\,
 \Big[\frac{1}{|\sin(\frac{\pi\, (x-y)}{2})|} + \frac{1}{|\sin(\frac{\pi\, (x+y)}{2})|}\Big]\, .
\end{equation}
Hence for $A\in ]0,2]$ and $\lambda\in ]0,\frac{2}{A+1}\wedge 1[$, we have for any $t>0$
\begin{equation}\label{B2(13)}
\sup_{x\in [0,1]} \int_0^1\! T_2(t,x,y)\, dy
\leq C\, e^{-ct}\, n^{-A\lambda}\, t^{-\frac{A+1}{2}\lambda}\, .
\end{equation}
In order to bound the $\|\;\|_{(\alpha)}$ norm of $T_2(t,x,.)$ for
$t\geq \gamma\, n^{-2}$, let $ N_1(n)=[\sqrt{n}]$, $
N_2(n)=[n/2]$ and  $ N_3(n)=n-1$. Then  $ T_{2}(t,x,y)\leq \sum_{i=1}^3T_{2,i}(t,x,y)$
where $
T_{2,1}(t,x,y)=\sum_{j=1}^{[\sqrt{n}]}|\Delta_{j}^n(t)| $ and for
$i=2,3$,  $
T_{2,i}(t,x,y)=\Big|\sum_{j=N_{i-1}(n)+1}^{N_i(n)}\!\!\Delta_{j}^n(t)
\varphi_{j}(x) \varphi_{j}(y)\Big|$. 
The inequalities (\ref{B2(10)}) with $A=2$ and (\ref{A4(1)}) with $\beta=0$  yield
$ 
\sup_{x,y\in [0,1]}T_{2,1}(t,x,y)\leq C\, \sum_{j=1}^{[\sqrt{n}]}n^{-1}e^{-cj^2t}
\leq C\, n^{-1}\,[1+ t^{-\frac{1}{2}}] \, e^{-ct}\, .
$ 
Furthermore, $ \sup\{  T_{2,1}(t,x,y)\,  ; \, (t,x,y)\in ]0,+\infty[\times [0,1]^2\}
\leq C \sqrt{n}$.
Hence both estimates yield
\begin{equation} \label{B2(15)}
\sup_{x\in [0,1]} \|T_{2,1}(t,x,.)\|^2_{(\alpha )} 
 \leq C\, e^{-ct}\, (1+t^{-1+\frac{\alpha}{3}})\, n^{-2+\alpha}. 
\end{equation}
For  $\lambda\in ]\alpha,1[$ and  $\mu\in ]0,1-\lambda[$, using (\ref{B2(12)})
 and (\ref{B2(11)}) with $A=0$, we have for $t\geq \gamma \, n^{-2}$, with the notations defined in
 (\ref{mathcalA}) and (\ref{mathcalB}):
 $ 
 \sup_{x\in [0,1]}
 T^{(1)}_{2,3}(t,x)
  \leq C n^{\alpha} e^{-ctn^2}[1+ (n t^{\frac{1}{2}})^{-(\lambda+\mu)}].
$ 
Similar computations for integrals over the sets ${\mathcal A}^{(i)}_n(x)$ for $i=2,3$ yield
$  \sup_{x\in [0,1]}  \|T_{2,3}(t,x,.)$\linebreak
$\times 1_{A_n^2(x)}(.)\|_{(\alpha)}^2 \leq C n^{\alpha} e^{-ctn^2}
.$
Furthermore, (\ref{B2(12)}) implies that
$ 
 \bar{T}^{(1)}_{2,3}(t,x) \leq C e^{-c n^2 t}
 \Big(\int_0^{2  n^{-1}} u^{-2} du\Big)$\linebreak 
$ \times  \Big(\int_0^{2 n^{-1}} v^{-\alpha} dv\Big) \leq  C e^{-c n^2 t} n^{\alpha}.$ 
For $(y,z)\in {\mathcal B}^{(2)}_n(x)$,  let
$I(y,z)\leq M(y,x)\leq S(y,z)$ denote the ordered values of $|x-y|$, $|y-z|$
 and $|x-z|$. Then 
$ 
 \bar{T}^{(2)}_{2,3}(t,x) \leq   C e^{-c n^2 t} \!\!
\int_{2 n^{-1}\leq I(y,z)\leq S(y,z)\leq 2}\! \Big(I(y,z)$\linebreak $\times  M(y,z) \Big)^{-1-\frac{\alpha}{2}}
  dy dz  
\leq  C e^{-c n^2 t}  n^{\alpha} .$ 
The previous inequalities on $T_{2,3}$ 
 yield that for $t\geq \gamma n^{-2}$, 
\begin{equation}\label{B2(16)}
\sup_{x\in [0,1]}  \|T_{2,3}(t,x,.)\|_{(\alpha)}^2 
\leq  C\, n^{\alpha}\, e^{-ctn^2}\,  .
\end{equation}
We now estimate $T_{2,2}$. Let $C_0>0$ be a "large" constant to be chosen later  on,  and 
 suppose that $t\geq C_0\, n^{-2}$.
For fixed $n$, $t$ and $j\in \big[[\sqrt{n}]+1,[n/2]\big]$, set
$\phi(j):=\Delta_{j}^n(t)$. Then
$\phi'(j)=2j\pi^2t\exp[-j^2\pi^2t]\; \big(
1-\psi(\frac{j\pi}{2n})\big)$, where for
$\frac{\pi}{2\sqrt{n}}\leq u \leq \frac{\pi}{4}$ one sets
$\psi(u):=\frac{\sin(2u)}{2u}\; \exp\big[\, 4n^2t\,
(u^2$\linebreak $-\sin^2u)\, \big]$. Hence, to apply Abel's
summation method, we have to compare $\psi(u)$ and 1.

\noindent  Using Taylor's expansion of the functions sine and exponential, we deduce that
 there exists a positive constant  $C_1$ such that if $\tilde{C}_1=(\frac{2\, C_1}{\pi})^2$,
for $C_0$ large enough,
the map $j\mapsto \phi(j)$ decreases on   $ \big[ [\sqrt{n}],  [n/2]\big]$ for $t\geq \frac{\tilde{C}_1}{n}$.
Let $t\in [\frac{C_0}{ n^2}, \frac{\tilde{C}_1}{n}]$. Then there exists a constant $C_2\in ]0,C_1[$ such that
$j\mapsto \phi(j)$ increases on $\big[ [\sqrt{n}],  [\frac{2\, C_2}{\pi\, \sqrt{t}}]\big]$ and decreases on
$\big[ [\frac{2\, C_1}{\pi\, \sqrt{t}}]+1,  [\frac{n}{2}]\big]$.
For $ t\in [\frac{C_0}{n^2}, \frac{\tilde{C}_1}{n}]$  $T_{2,2}(t,x,y)\leq \sum_{i=1}^2 T_{2,2,i}(t,x,y)$, where
one set $B_i=\frac{2C_i}{\pi}$ and
\begin{align*}
T_{2,2,1}(t,x,y)=&\Big|\sum_{j=[\sqrt{n}]}^{[B_2 / \sqrt{t}]} \Delta_j^n(t)\, \varphi_j(x)\,
\varphi_j(y)\Big|
+ \Big|\sum_{j=[B_1/ \sqrt{t}]}^{[n/2]} \Delta_j^n(t)\, \varphi_j(x)\,
  \varphi_j(y)\Big|\, ,\\
 T_{2,2,2}(t,x,y)=&\Big|\sum_{j=[B_2 / \sqrt{t}]}^{[B_1/ \sqrt{t}]} \Delta_j^n(t)\, \varphi_j(x)\,
\varphi_j(y)\Big|\, .
\end{align*}
There exists  a constant $C>0$ such that if $\gamma > C_0$,  for every  $t \geq \frac{\gamma}{n^2}$
\begin{equation}\label{trivial}
\sup_{(x,y)\in [0,1]^2} T_{2,2,1}(t,x,y)\leq C\, n\, e^{-ct}.
\end{equation}
For $t\geq \frac{\tilde{C}_1}{n}$, let $T_{2,2,1}(t,x,y)=T_{2,2}(t,x,y)$.
 Using (\ref{B2(12)})  we deduce that for $t\geq \frac{C_0}{n^2}$ and $\beta\in [0,1]$
\begin{align}\label{majoT221}
T_{2,2,1}(t,x,y) \leq & C \Big[\frac{1}{|\sin(\frac{\pi (x-y)}{2})|}+
 \frac{1}{|\sin(\frac{\pi (x+y)}{2})|}\Big]
\Big[ \Delta_{[\sqrt{n}]+1}^n(t) +\sum_{i=1}^2 
\Delta_{[\frac{{B_i}}{\sqrt{t}}]}^n(t) 
1_{\{C_0n^{-2}\leq t\leq \tilde{C}_1 n^{-1}\}}   \Big] \nonumber \\
 \leq & C  \Big[\frac{1}{|x-y|}+ \frac{1}{x+y}+ \frac{1}{2-x-y}\Big]
\Big[  \frac{e^{-ctn}}{n} + n^{-2\beta} t^{-\beta} 1_{C_0 n^{-2}
\leq t\leq \tilde{C}_1n{-1} \}} \Big]  .
\end{align}
For  $t\in [\frac{C_0}{ n^2}, \frac{\tilde{C}_1}{n}]$, it  remains to  
bound directly the sum $T_{2,2,2}(t,x,y)$. The inequality (\ref{B2(10)})
implies that for $B_2\, t^{-\frac{1}{2}}\leq j\leq B_1\, t^{-\frac{1}{2}}$,
$\Delta_j^n(t)\leq C n^{-2A}\, t^{-A}\, e^{-c\, t\, j^2}$ for any $A\in [0,1]$. Therefore,
the inequality (\ref{A4(1)}) implies that for  $A\in [0,1]$,
\begin{equation}\label{majoTun222}
\sup_{x\in [0,1]} T_{2,2,2}(t,x,y)\leq C\, n^{-2A}\, t^{-(A+\frac{1}{2})}\, .
\end{equation}
Finally, if $C_0$ is large enough, for $t\in [\frac{C_0}{n^2}, \frac{\tilde{C}_1}{n}]$  
the function $\psi$ is increasing on the interval $[\frac{C_2}{ n \sqrt{t}}\, ,
\frac{ C_1}{ n \sqrt{t}}]$
and  $\sup\{\phi(u)\, :\, u\in [\frac{C_2}{n\sqrt{t}},\frac{C_1}{n\sqrt{t}}]\}
\leq C\, n^{-2}\, t^{-1}$. Hence Abel's summation method implies that for $t\in [\frac{C_0}{n^2}, \frac{\tilde{C}_1}{n}]$,
\begin{equation}\label{majoT222trois}
\sup_{x\in [0,1]} T_{2,2,2}(t,x,y)\leq C\,  n^{-2}\, t^{-1} \Big[\, \frac{1}{|\sin(\pi\, \frac{ x-y}{2})|}+
\frac{1}{|\sin(\pi \, \frac{x+y}{2})|}\, \Big]\, .
\end{equation}
 The inequalities (\ref{majoT221}) applied with $\beta=\frac{1}{2}$ and $\beta=1$ respectively and (\ref{trivial})
imply that for $\lambda\in ]0,\alpha[$ and  $\mu\in ]0,1-\lambda[$, $\nu=\lambda+\mu$,
there exists a constant $C>0$ such that for every
 $t\geq \gamma n^{-2}$:
$  
 \sup_{x\in [0,1]}  
T^{(1)}_{2,2,1}(t,x)\leq 
 C \Big[n^{-1+\alpha} e^{-ctn} + 1_{\{C_0 n^{-2}\leq t\leq \tilde{C}_1 n^{-1}\}}
n^{-4+\alpha+2\nu} t^{-\nu}\Big]  .
$
Similar computations for the integrals over the sets ${\mathcal A}^{(i)}_n(x)$, $i=2,3$ imply that
the same upper estimates   hold for $T^{(i)}_{2,2,1}(t,x)$.
Furthermore, (\ref{majoT221}) with $\beta\in ]\frac{1}{2},1[$ yields
$
 \sum_{j=1}^2 \sup_{x\in [0,1]}  \bar{T}^{(j)}_{2,2,1}(t,x)
\leq \, C\, \big( n^{-2+\alpha}$ \linebreak[4] $\times  e^{-ctn} + n^{-4\beta +\alpha}\, t^{-2\beta}\,
1_{\{C_0\, n^{-2}\leq t\leq \tilde{C}_1\, n^{-1}\}} \big) \, .
$ 
The inequalities 
on ${T}^{(i)}_{2,2,1}$ and $\bar{T}^{(j)}_{2,2,1}$
 yield that for  $\nu\in ]\alpha,1[$ and $\beta\in ]\frac{1}{2},1[$,
\begin{equation}\label{B2(18-14)}
\sup_{x\in [0,1]} \|T_{2,2,1}(t,x,.)\|_{(\alpha)}^2 \leq C \big[ n^{-1+\alpha} e^{-ctn}+
1_{\{\frac{C_0}{n^2}\leq t\leq \frac{\tilde{C}_1}{n}\}} \big(  t^{-\nu}n^{-4+\alpha+2\nu}
+ t^{-2\beta}n^{-4\beta +\alpha}\big) \big] .
\end{equation}
For $t\in [\frac{C_0}{ n^2}, \frac{\tilde{C}_1}{n}]$  
  (\ref{majoTun222}) and (\ref{B2(12)}) yield
that  for  $A\in [0,1]$,  $\lambda\in ]0,\alpha[$,  $\mu\in ]0,1[$ and $\nu=\lambda + \mu$,
$  
\sup_{x\in [0,1]}\sum_{i=1}^3 T_{2,2,2}^{(i)}(t,x)
 \leq C n^{-(2A+1)\nu+\alpha} t^{-(A+\frac{1}{2})\nu}.
$ 
Proceeding as for the estimates of $\bar{T}^{(i)}_{2,3}$ we deduce that for
 $A =1 $ and $\lambda\in ]\frac{1}{3},\frac{1}{2}[$
(resp. for  $0<\tilde{\lambda}<\frac{\alpha}{2}$),
$ \sup_{x\in [0,1]} \bar{T}^{(1)}_{2,2,2}(t,x)
\leq C\, n^{-6\lambda+\alpha}\, t^{-3 \lambda}$ and  
$ \sup_{x\in [0,1]} \bar{T}_{2,2,2}^{(2)}(t,x)
\leq  C\, n^{-4}\, t^{-2-\tilde{\lambda}} \,   \Big(\int_{2\, n^{-1}}^2
u^{-1+\tilde{\lambda}-\frac{\alpha}{2}}\, du\Big)^2
\leq C\, n^{-4+\alpha}\, t^{-2}$.

\noindent The upper estimates of ${T}_{2,2,2}^{(i)}(t,x)$,  $ \bar{T}_{2,2,2}^{(j)}(t,x)$
 and (\ref{B2(18-14)})  imply that for
$\gamma$ large enough, $\mu\in ]\frac{1}{3},\frac{1}{2}[$  there exists  $C>0$ such that for every
$t\geq  \gamma n^{-2}$:
\begin{equation}\label{B2(18)}
\sup_{x\in [0,1]}\|T_{2,2}\|_{(\alpha)}^2 \leq Cn^{\alpha}  \Big[ e^{-ctn} n^{-1} +
1_{\{\frac{C_0}{n^{2}}\leq t\leq \frac{\tilde{C_1}}{ n}\}}
\Big( n^{-4+2\nu} t^{-\nu} +  t^{-2}n^{-4 } +  t^{-\lambda} n^{-2\lambda } \Big)\Big] .
\end{equation}
Thus (\ref{B2(15)}), (\ref{B2(16)}) and   (\ref{B2(18)})  yield for
  $\lambda \in ]1,\frac{3}{2}[$,  $\nu\in ]\alpha,1[$ and $t\geq \gamma n^{-2} $ with $\gamma$
large enough, 
\begin{align}\label{B2(19)}
 \sup_{x\in [0,1]} \| T_{2}(t,x,.)\|_{(\alpha)}^2  \leq &\,  Cn^\alpha \Big[
e^{-ctn^2}+ n^{-2} e^{-ct} (1+t^{-1+\frac{\alpha}{3}}) +
n^{-1} e^{-ctn}\nonumber \\
&\:  + 1_{\{\frac{C_0}{n^2}\leq t\leq \frac{\tilde{C_1}}{ n}\}}
\Big( n^{-4+2\nu} t^{-\nu}+ n^{-4  } t^{-2} + n^{-2\lambda } t^{-\lambda}\Big)\Big]  .
\end{align}
\noindent {\bf Estimates of $T_{3}$}
Using (\ref{A4(1)}) we deduce that
$T_{3}(t,x,y)\leq  C\, n^{-1}\,
 \big[1+t^{-1}\big]\, e^{-ct}\, .$
 Furthermore,  set
$A(l):=[0,1]\cap\big(\big[\frac{l-1}{n},\frac{l+2}{n}\big]\cup\big[0,\frac{(2-l)^+}{n}\big]
\cup \big[(2-\frac{l+2}{n})\wedge 1,1\big]\big)$. Then $dx(A(l))\leq C\, n^{-1}$.
The study of the monotonicity of the function $H$ defined by
$ H(z)=z\, \exp\left[-4n^2t\sin^2\left(\frac{z\pi}{2n}\right)\right]$ and
 Abel's summation method yield
for large enough $\gamma $,  $t\geq \gamma\, n^{-2}$, $0<\lambda<1$ and $y\in A(l)$,
$\sup_{x\in [\frac{l}{n},\frac{l+1}{n}]} T_3(t,x,y)
\leq C\,  (1+t^{-\frac{1+\lambda}{2}})\, n^{-\lambda}\, e^{-ct}\, ,
$ 
while for $y\not\in A(l)$,
$\sup_{x\in [\frac{l}{n},\frac{l+1}{n}]} T_3(t,x,y)\leq \frac{C}{n}\,  (1+t^{-\frac{1}{2}})\, e^{-ct}\,
\big[\; \big|y-\frac{2l+1}{2n}\big|^{-1}+ \big|y+\frac{2l+1}{2n}\big|^{-1}+ \big|2n-y-\frac{2l+1}{2n}\big|^{-1} \, \big]\, .$
Then, using the partition $\{A(l),A(l)^c\}$ and the three previous inequalities
we deduce that for  $\bar{\lambda}\in ]0,1[$ and  $t\geq \gamma\, n^{-2}$ one has
\begin{equation}\label{B2(31)}
\sup_{x\in [0,1]} T_3(t,x,y)\leq 
C e^{-ct} (1+t^{-\frac{1+\bar{\lambda}}{2}}) n^{-1}.
\end{equation}
Similarly, for $t\geq \gamma\, n^{-2}$ and $\lambda\in ]0,1[$ 
there exists a constant $C>0$
such that for every  $l\in \{0, \cdots,n-1\}$ and $x\in [\frac{l}{n},\frac{l+1}{n}]$,
$ 
\|T_3(t,x,.) 1_{A(l)}(.)\|_{(\alpha)}^2\leq  
C e^{-ct} (1+t^{-1-\lambda}) n^{-2+\alpha-2\lambda}.
$ 
Furthermore, when  $t\geq \gamma\, n^{-2}$ 
 separate  estimates
in the cases $y,z\not\in A(l)$ and either $|y-z|\leq n^{-1}$ or $|y-z|\geq n^{-1}$ yield that
given $\nu \in ]0,\frac{\alpha}{2}[$,  there exists a constant $C>0$ such that
for every $l\in \{0, \cdots,n-1\}$ and $x\in [\frac{l}{n},\frac{l+1}{n}]$, one has 
$  
\ddd \|T_3(t,x,.) 1_{A(l)^c}(.)\|_{(\alpha)}^2
 \leq C n^{-2-\nu+\alpha} e^{-ct} (1+t^{-\frac{2+\nu}{2}}).
$  
These inequalities 
 imply that for $t\geq \gamma\, n^{-2}$ and  $\lambda \in ]0,\frac{\alpha}{4}[$, 
\begin{equation}\label{B2(32)}
\sup_{x\in[0,1]} \|T_3(t,x,.)\|_{(\alpha)}^2\leq C\, e^{-ct}\, n^{-2+\alpha-2\lambda}\, (1+t^{-1-\lambda}) \, .
\end{equation}
\noindent {\bf Estimates of  $T_{4}$}  We  suppose that $ x=\frac{l}{n}$, $1\leq l\leq n-1$.
The general case is  easily deduced by linear interpolation. For $\frac{k}{n}\leq y\leq \frac{k+1}{n}$,
$0\leq k\leq n-1$, one has $ k_{n}(y)=\frac{k}{n}$ and using  (\ref{A4(1)}), we deduce
$ T_{4}(t,x,y) 
\leq C\, n^{-1}\, [\, t^{-1}+1\,]\, e^{-ct}\, .
$ 
Let $B(l):=\{u\in [0,1]; \left|\frac{l}{n}-u\right|\leq\frac{1}{n}\mbox { or }\frac{l}{n}-u\leq\frac{1}{n}
\mbox { or }2-\frac{l}{n}-u\leq\frac{1}{n}\}$;  as usual, $dx(B(l))\leq cn^{-1}$.
  Let then $C^{1}(l):=\{y\in [0,1];\exists u\in B(l)\cap[k_{n}(y),y]\}$ and for $i=1,2$, let 
$\tilde{C}^{i}(l):=\{z\in [0,1];\exists y\in C^{1}(l),|y-z|\leq\frac{i}{n}\}$. Then $dx(\tilde{C}^{i}(l))\leq C\, n^{-1}$
and for $y\not\in\tilde{C}^{1}(l)$, one has $|y-x|\wedge(y+x)\wedge(2-x-y)\geq n^{-1}$.
Computations 
similar to that made to estimate $T_3$ yield for $\lambda\in ]0,1[$
the existence of a constant $C>0$ such that
for every $l\in \{0, \cdots, n\}$, and $y\in \tilde{C}^2(l)$,
 $T_4(t,l/n,y) 
\leq  C n^{-\lambda}\big(1+t^{-\frac{1+\lambda}{2}}\big) e^{-ct}$ 
and for $y\in \tilde{C}^2(l)^c$, $T_4(t,l/n,y) 
 \leq  C n^{-1} (1+t^{-\frac{1}{2}}) e^{-ct} \big[\, \big|y-\frac{l}{n}\big|^{-1}
+ \big(y+\frac{l}{n}\big)^{-1} + \big|2-y-\frac{l}{n}\big|^{-1}\big] .$
An argument similar to that proving (\ref{B2(31)}) and (\ref{B2(32)}) implies that for  
$\lambda\in ]0,1[$ and $\nu\in ]0,\frac{\alpha}{4}[$, there exists a constant $C>0$ such that
for any $t\geq \gamma \, n^{-2}$
\begin{align}
\sup_{x\in [0,1]}T_{4}(t,x,y)\leq& C  e^{-ct} (1+t^{-\frac{1+\lambda}{2}}) n^{-1}, \label{B2(34)}
\\
\label{B2(35)}
\sup_{x\in [0,1]}\|T_{4}(t,x,\cdot)\|_{(\alpha)}^2 \leq&   n^{-2+\alpha} e^{-ct}
 (1+t^{-(1+\nu)}) n^{-2\nu}  .
\end{align}
\noindent {\bf Estimate of $T_1(t,x,.)$.} Using (\ref{A4(1)}) with $\beta=0$ and $J_{0}=n$, we have
\begin{equation}\label{numero1}
\sup_{x,y\in [0,1]} T_{1}(t,x,y) \leq C\, \sum_{j\geq n}e^{-ctj^2}\leq C\, e^{-ctn^2}\, [1+ t^{-\frac{1}{2}}].
\end{equation}
On the other hand, since $j\rightarrow e^{-j^2\pi^2t}$ decreases,  Abel's summation method yields
$ 
T_1(t,x,y)\leq C\, e^{-cn^2t}\, \big[\frac{1}{|\sin(\pi \frac{x-y}{2})|} +
 \frac{1}{|\sin(\pi \frac{x+y}{2})|} \big]\, .
$ 
Thus, for $\lambda\in [0,1[$ and $t>0$ we have
\begin{equation}\label{B2(6)}
\sup_{x\in [0,1]} \int_0^1T_1(t,x,y)\, dy 
\leq C\, e^{-cn^2t}\, [1+t^{-\frac{\lambda}{2}}]\, .
\end{equation}
Given  $\lambda\in ]\alpha,1[$, $\mu\in 
]0,1[$ and $t\geq \gamma n^{-2}$, computations similar to the previous ones
using 
the partition $\{ A_n^2(x), A_n^2(x)^c\}$  yield that for 
 $t\geq \gamma\, n^{-2}$, we have
\begin{equation}\label{B2(8)}
\sup_{x\in [0,1]} \| T_1(t,x,.)\|_{(\alpha)}^2 \leq C\, n^\alpha\, e^{-cn^2t} \, .
\end{equation}
The inequalities (\ref{B2(6)}), (\ref{B2(13)}) with $A=2$, $\lambda=\frac{1}{2}$,
and $(T_3+T_4)(t,x,y)\leq Cn^{-1}(1+t^{-1}) e^{-ct}$\ 
 imply the existence of  $c,C>0$, $\lambda\in ]0,1[$
such that for   $t>0$,
\begin{equation}\label{B2(36bis)}
 \sup_{x\in [0,1]} \|G(t,x,.)-G^n(t,x,.)\|_1 \leq C\, n^{-1}\, \Big[ \big(
1+t^{-1}\big)\, e^{-ct}
+t^{-\frac{\lambda}{2}}\, e^{-ctn^2}\Big]\, ,
\end{equation}
which implies (\ref{B2(0)}) for $d=1$ with $\mu=1$.
The inequalities (\ref{B2(6)}), (\ref{B2(13)}) with $A=2$ and $\lambda=\frac{1}{2}$, (\ref{B2(31)})
 and (\ref{B2(34)}) with $\bar{\lambda}=\frac{1}{2}$ imply that for some $\mu \in ]0,\frac{1}{2}[$ there exists
 a constant $C>0$ such that for every $t\geq \gamma\, n^{-2}$
with  $\gamma >0$ large enough, one has
\begin{equation}\label{B2(36)}
\sup_{x\in[0,1]} \|G(t,x,.)-G^n(t,x,.)\|_1 \leq C\,  \Big[\, n^{-1}\, \Big(1+t^{-\frac{3}{4}}\Big) +
e^{-c\, t\, n^2} \, (1+t^{-\mu})\, \Big]\, .
\end{equation}
On the other hand, the inequalities (\ref{B2(8)}), (\ref{B2(19)}), (\ref{B2(32)}) and (\ref{B2(35)}) imply that for
$\nu\in ]0,\frac{\alpha}{4}[$, $\lambda \in ]1,\frac{3}{2}[$, $\mu\in ]\alpha, 1[$, there exist  positive constants 
$C$ and $\tilde{C}_1$  such that   for  $t\geq \gamma \, n^{-2}$ with  $\gamma >0$ large enough, one has
\begin{align}\label{B2(37)}
&\sup_{x\in[0,1]} \|G(t,x,.) -G^n(t,x,.)\|_{(\alpha)}^2 \leq C n^{\alpha}  \Big[ e^{-ctn^2}
+ \frac{e^{-ctn}}{n}
+\frac{e^{-ct}}{n^2}  \Big(1+ \frac{t^{-(1+\nu)}}{n^{2\nu}} + t^{-1+\frac{\alpha}{3}}\Big) \nonumber\\
&\qquad 
 + \Big( n^{-4+2\mu} t^{-\mu}+n^{-4} t^{-2}+n^{-2\lambda} t^{-\lambda}\Big)
 1_{\{\gamma n^{-2}\leq t\leq \tilde{C}_1 n^{-1}\}} \Big] ,
\end{align}
which proves 
(\ref{B2(0)})-(\ref{B2(2)})  for $d=1$.

\noindent {\bf The general case}  We extend the inequalities  (\ref{B2(36bis)}),
 (\ref{B2(36)}) and (\ref{B2(37)}) to any dimension
$d$. We use the fact that for any $d\geq 2$, we have $|G_d(t,x,y)-(G_d)^n(t,x,y)|
 \leq \sum_{i=1}^d  \Pi_i$, where 
\begin{align}\label{B2(38)}
&
\Pi_i= \Big(\prod_{j=1}^{i-1} |G(t,x_j,y_j)|\Big) 
 |G(t,x_i,y_i)-G^n(t,x_i,y_i)|  \Big(\prod_{j=i+1}^d |G^n(t,x_j,y_j)|\Big) .
\end{align}
Hence, the inequalities (\ref{A1(0)}), (\ref{B1(1)}),
(\ref{B2(36bis)}), (\ref{B2(36)}) and (\ref{B2(38)}) imply (\ref{B2(0)}) with
some $\mu>0$, and
that for any  $\lambda\in ]0,1[$   
and any $\nu\in ]0,1/4[$,
there exists  $C>0$ such that  for   $t\geq \gamma\, n^{-2}$,
\begin{equation}\label{B2(39)}
\sup_{x\in Q}\|G_d(t,x,.)-(G_d)^n(t,x,.)\|_1\leq C
 t^{-\nu} \big[ \big( 1+t^{-\frac{3}{4}}\big)n^{-1} +
e^{-c t n^2}  (1+t^{-\nu}) \big].
\end{equation}
Using (\ref{B2(3)}) and integrating (\ref{B2(39)}) with respect to $t$ on $[\gamma\, n^{-2},+\infty[$
we obtain (\ref{B2(1)}).
Finally, for $1\leq k\leq d-1$ set $\alpha_k=\alpha\, 2^{-k}$ and set  $\alpha_d=\alpha_{d-1}$; then using
(\ref{B2(38)}), (\ref{majof}), (\ref{A1(1)}), (\ref{B1(2)}) and (\ref{B2(37)}),  we deduce that for $\alpha\in ]0,2[$,
$C_1>0$, $\lambda\in ]\alpha ,1[$, $\mu\in ]1,\frac{3}{2}[$, $\nu\in ]0,\frac{\alpha_d}{4}[$
   and for  $t\geq \gamma \, n^{-2}$ for $\gamma >0$ large enough 
\begin{align}\label{B2(40)}
& \sup_{x\in Q} \| G_d(t,x,.-(G_d)^n(t,x,.) \|_{(\alpha)} \leq C\,\sum_{i=1}^d t^{-\frac{1}{2}\, \sum_{j=1}^{i-1}\alpha_j}
\nonumber \\
& \qquad\qquad \times \|G(t,x_i,.)-G^n(t,x_i,.) \|_{(\alpha_i)} \, n^{\sum_{j=i+1}^d \alpha_j}\, \nonumber\\
&\quad \leq C\, n^{\alpha}   \Big[\, e^{-ctn^2}
+ n^{-1}\, e^{-ctn} + n^{-2}\, e^{-ct}\, \big( t^{-(1+\nu)}\, n^{-2\nu}+ t^{-(1+\lambda)}\, n^{-2\lambda}+1\big) \nonumber \\
&
\qquad \qquad  + \big( n^{-4}\, t^{-2}+n^{-4+2\nu}\, t^{-2+\nu}+n^{-2\mu}\, t^{-\mu}\big)\,
 1_{\{\gamma n^{-2}\leq t\leq \tilde{C}_1 n^{-1}\}}\, \Big]\, .
\end{align}
Integrating  on  $[\frac{\gamma}{ n^2},+\infty[$ and using (\ref{B2(4)}),
we deduce (\ref{B2(2)}).
\qquad $\Box$

\noindent We now estimate the norm of the difference $(G_d)^n$ and $(G_d)^{n,m}$.
\begin{lemme}
\label{B3}
Given any $T>0$ and  $\nu >0$ there exists $C>0$ such that
\begin{align}
&\sup_{x\in Q} \int_{0}^T\int_Q \big[\,  |(G_d)^n(t,x,y)-(G_d)^{n,m}(t+T\, m^{-1},x,y)| \nonumber \\
&\qquad + |(G_d)^n(t,x,y)-(G_d)^n_m(t+T\, m^{-1},x,y)|\, \big]\, dydt \leq C\, m^{-1+\nu}\, ,
\label{B3(3)}\\
 \lefteqn{\sup_{x\in Q} \int_{0}^T\big[\, \| (G_d)^n(t,x,\cdot)-(G_d)^{n,m}(t+T\, m^{-1},x,\cdot) \|_{(\alpha)}^2}\nonumber\\
&\qquad + \|(G_d)^n(t,x,.)-(G_d)^n_m(t,x,.)\|_{(\alpha )}^2\,\big]\,  dt\leq
C\,  m^{-1+\frac{\alpha}{2}}\, . \label{B3(4)}
\end{align}
\end{lemme}
{\bf Proof : } We only prove these inequalities for $G^n-G^{n,m}$ and
we at first suppose that $d=1$. Let $\bar{\mathcal
G}^{n,m}=G^n-\tilde{\mathcal G}^{n,m}$ and
$\bar{G}^{n,m}=G^{n,m}-\tilde{G}^{n,m}$ where $\tilde{G}^{n,m}$ is
defined
by (\ref{tildeG}), 
and $\tilde{\mathcal G}^{n,m}$ is defined by
\begin{equation}\label{tildeG2}
\tilde{\mathcal G}^{n,m} (t,x,y)=\sum_{j=1}^{(n\wedge \sqrt{m})-1}
e^{\lambda_j^n t}\, \varphi_j^n(x) \varphi_j(\kappa_n(y))\, .
\end{equation}
Then (\ref{B1(6-1)}) and (\ref{B20}) provide upper estimates of
the norms of $\bar{G}^{n,m}$. Similar computations prove that the
same upper estimates hold for the norms of $\bar{\mathcal
G}^{n,m}$, i.e.,
 for $\lambda\in ]0,\frac{1}{2}[$ and $\beta\in ]\alpha,1[$,
\begin{align}\label{B3(7)}
&\sup_{x\in [0,1]} \|\bar{\mathcal G}^{n,m}(t,x,.)\|_1\leq C
 (1+t^{-\lambda}) e^{-ctm}\, , \;  
\sup_{x\in [0,1]} \| \bar{\mathcal G}^{n,m}(t,x,.) \|_{(\alpha)}^2
\leq C m^{\frac{\alpha-\beta}{2}}
t^{-\frac{\beta}{2}}e^{-ctm}.
\end{align}
Let $\tilde{c}$ be a positive constant to be fixed later on; for
$t\leq \tilde{c}\, T\, m^{-1}$ we estimate separately the norms of
$\tilde{G}^{n,m}(t,x,.)$ and $\tilde{\mathcal G}^{n,m}(t,x,.)$.
The inequalities (\ref{B1(4-1)}) and (\ref{B1(4-2)}) provide  the
estimates of $\tilde{G}^{n,m}$.  For $\tilde{\mathcal G}^{n,m}$, we
proceed in a similar way. Indeed, $j\rightarrow \exp(\lambda_j^n\,
t)$ is decreasing,  $\exp(\lambda_j^n\, t)\leq e^{-ctj^2}$ for
$c>0$, and $|\tilde{\mathcal G}^{n,m}(t,x,y)|\leq C \,(
n\wedge\sqrt{m})$.
Hence the arguments used in the proof of Lemma \ref{B1mn} yield
that for any $\tilde{c}>0$ there exists a constant $C>0$ such that
for $t\in ]0,  \frac{\tilde{c} T}{ m}]$,
\begin{align}\label{B3(8-1)}
&\sup_{x\in [0,1]}\|\tilde{\mathcal G}^{n,m}(t,x,.)-\tilde{G}^{n,m}(t,x,.)\|_1 \leq C\, (1+t^{-\lambda})\, ,\\
\label{B3(8-2)} &\sup_{x\in [0,1]}\|\, \tilde{\mathcal
G}^{n,m}(t,x,.)-\tilde{G}^{n,m}(t,x,.) \|_{(\alpha)}^2 \leq C\,
(n\wedge \sqrt{m})^\alpha\, .
\end{align}
Furthermore, if $t\in [\tilde{c}\, T\, m^{-1},T]$,
$|\tilde{\mathcal G}^{n,m}(t,x,y)-\tilde{G}^{n,m}(t,x,y)|\leq
\tilde{T}(t,x,y)= \tilde{T}_{1}(t,x,y)+\tilde{T}_{2}(t,x,y)$,
where ${\displaystyle \tilde{T}_{i}(t,x,y) =
\big|\!\!\! \sum_{j=1}^{(n\wedge\sqrt{m})-1} A^i_{n,m}(t)
\varphi_{j}^n(x)\varphi_{j}(\kappa_{n}(y))\big|}$, 
$A^1_{n,m}(t) =   \exp\left(
\frac{\left(\left[\frac{mt}{T}\right]+1\right) \,
\lambda_{j}^n\,T}{m}\right)$\linebreak $
 -\left(1-\lambda_{j}^n\, \frac{T}{m}\right)^{-\left(\left[\frac{mt}{T}\right]+1\right)}
$ and $A^2_{n,m}(t,x)= \exp\left(\lambda_{j}^nt\right)-
\exp\left(\frac{\left(\left[\frac{mt}{T}\right]+1\right)\,
 \lambda_{j}^n\, T}{m}\right)$.
 Using Abel's summation method, we have for $i=1,2$,  $t\geq\frac{\tilde{c}}{m}$ for $\tilde{c}$ large enough,
 $ x=\frac{l}{n}$ and $\kappa_n(y)=\frac{k}{n}$,
$ |\tilde{T}_{i}(t,x,y)|
\leq C \frac{T}{mt} \Big[
\frac{1}{|x-\kappa_n(y)|}+\frac{1}{x+\kappa_n(y)}+
\frac{1}{2-x-\kappa_n(y)}\Big].$
Furthermore,
$ 
|\tilde{T}(t,x,y)|\leq 
C (n\wedge \sqrt{m}).
$ 
These inequalities  
 yield  that 
 for $t\geq \frac{\tilde{c}T}{m}$ for large enough $\tilde{c}$  and $\lambda\in ]0,1[$,
\begin{equation}\label{B3(12)}
 \|\tilde{T}(t,x,.)\|_1\leq C t^{-1+\lambda}
m^{-1+\frac{3\lambda}{2}}.
\end{equation}
For $\mu \in ]\alpha,1[$, $\nu\in ]0,1-\mu[$ and $\beta
=\mu+\nu\in ]\alpha,1[$, using the sets
$\mathcal{A}_{n\wedge\sqrt{m}}^{(i)}(x)$ for $i\leq 3$ and
$\mathcal{B}_{n\wedge \sqrt{m}}^{(j)}(x)$ for $j=1,2$ and the fact
that $\frac{1}{n}\leq \frac{1}{n\wedge\sqrt{m}}$, 
we deduce that given $\tilde{c}$ large enough, there exists  constants
$c,C>0$ such that for every $t\in [\frac{\tilde{c}T}{m},T]$:
\begin{align}\label{B3(13)}
&\| \tilde{T}(t,x,.) \|^2_{(\alpha)} 
\leq C\, (n\wedge\sqrt{m})^{\alpha}\,m^{-2}\, t^{-2} \leq C\, m^{-2+\frac{\alpha}{2}}\, t^{-2}\, .
\end{align}
For $d=1$, the inequalities (\ref{B1(6-1)}), (\ref{B3(7)}),
(\ref{B3(8-1)}) and (\ref{B3(12)}) imply the existence of $\lambda
\in ]0,\frac{1}{2}[$ and positive constants $c,C$ such that for
any $t\in ]0,T]$:
\begin{equation}\label{B3(1)}
\sup_{x\in Q}\|\big( (G_d)^n-(G_d)^{n,m}\big)(t,x,.)\|_1 \leq C\, \Big[ \big(1+t^{-\lambda}\big)\, e^{-ctm}
+ t^{-1+\lambda}\, m^{-1+\frac{3\lambda}{2}}\Big]\, ,
\end{equation}
while the inequalities (\ref{B20}), (\ref{B3(7)}), (\ref{B3(8-2)}) and (\ref{B3(13)}) yield
the existence of $\beta \in ]0, \alpha \wedge d[$ and positive constants $\tilde{c}, c $ and $C$ such that for
every $t\in ]0,T]$:
\begin{equation}\label{B3(2)}
\sup_{x\in Q}\|\big( (G_d)^n-(G_d)^{n,m}\big)(t,x,.) \|_\alpha^2 \leq C m^{\frac{\alpha}{2}}
\Big[ \big( 1+ (tm)^{-\frac{\beta}{2}}\big){e^{-ctm}} +  1_{[\tilde{c} T m^{-1},T]}(t){m^{-2} t^{-2}}\Big].
\end{equation}
Let $\alpha_k=\alpha \,2^{_k}$ for $1\leq k\leq d-1$ and $\alpha_d=\alpha_{d-1}$.
The inequalities (\ref{B3(1)}) for $d=1$, (\ref{B1(1)}) and (\ref{B1(3-1)}) yield (\ref{B3(1)}) for any $d$, while
 (\ref{B3(2)}) for $d=1$, (\ref{B1(2)}) and (\ref{B1(3-2)}) yield (\ref{B3(2)}) for any $d$.
Integrating with respect to $t$ we deduce the inequalities (\ref{B3(3)}) and (\ref{B3(4)}).
\qquad$\Box$

\section{Some numerical results}\label{numerique}
In order to study the influence of the correlation coefficient $\alpha$ of the Gaussian noise
on the speed of convergence, we have implemented in C the implicit discretization scheme
$u^{n,m}$ in the case of homogeneous boundary conditions in dimension $d=1$ for the equation
(\ref{implicit2}).

To check the influence of the time mesh, we have fixed the space mesh  $n^{-1}$ with $n=500$
and taken the smallest time mesh $m_0^{-1}$ with $m_0=20736$. Using one trajectory of the noise $F$,
we have approximated by the Monte-Carlo method $e(m_i)=\E(|u^{n,m_0}(1,.5)-u^{n,m_i}(1,.5)|^2)$
and $\hat{e}(m_i)=\sup_{x\in [0,1]}\E(|u^{n,m_0}(1,x)-u^{n,m_i}(1,x)|^2)$ for 13 divisors $m_i$
of $m_0$, ranging from $m_1=854$ to $m_{13}=144$. These simulations have been done for various
values of $\alpha$, including the case of the space-time white noise. 
 Assuming that $u^{n,m_0}$ is close to $u$, according to (\ref{vitesnc3}) and (\ref{convn-nm}),
these errors should behave like $C\,[ m_i^{-(1-\frac{\alpha}{2})}+n^{-(2-\alpha)}]\sim m_i^{-(1-\frac{\alpha}{2})}$
for this choice of $n$ and $m_i$. Thus, we have computed the linear regression coefficients $c(t)$ and $d(t)$
(resp. $\hat{c}(t)$ and $\hat{d}(t)$) of $\ln(e(m_i))$ (resp. of $\ln(\hat{e}(m_i))$), i.e., of the approximation
of $\ln(e(m_i))$ by $c(t) \ln(m_i)+d(t)$ as well as the corresponding standard deviation $sd$ (resp. $\hat{sd}$)
for $K=3200$ Monte-Carlo iterations in the case $\sigma(x)=0.2\, x +1$ and $b(x)=x+2$. 

The study of the influence of the space mesh is done in a similar way; we fix the time mesh $m^{-1}$  with $m=32000$ and let
  the smallest space mesh $n_0=432$. Again for various divisors of $n_0$, using one trajectory of the noise $F$ we have
  approximated $\varepsilon(n_i)= \E(|u^{n_0,m}(1,.5)-u^{n_i,m}(1,.5)|^2)$ and $\tilde{\varepsilon}(n_i)=
  \E(|u^{n_0,m}(1,.5)-u^{n_i,m}(1,.5)|^2)$ for the 7 divisors $n_i$ of $n_0$ ranging from $72$ to $12$.
  Assuming that $u^{n_0,m}$ is close to $u$, according to (\ref{vitesnc3}) and (\ref{convn-nm}),
  these errors should behave like $C\,[ m^{-(1-\frac{\alpha}{2})}+n_i^{-(2-\alpha)}]\sim n_i^{-(2-\alpha)}$
  for this choice of $n_i$ and $m$. Thus, we have computed the linear regression coefficients $\gamma(x)$ and $\delta(x)$
  (resp. $\hat{\gamma}(t)$ and $\hat{\delta}(t)$) of $\ln(\varepsilon(n_i))$ (resp. of $\ln(\hat{\varepsilon}(n_i))$),
  i.e., of the approximation of $\ln(\varepsilon(n_i))$ by $\gamma(x) \ln(n_i)+\delta(x)$ as well as the corresponding
  standard deviation $SD$ (resp. $\hat{SD}$) for $K=3200$ iterations in the case $\sigma(x)=1$ and $b(x)=2x+3$.
 Both sets of results  are summarized as follows.
{\small \begin{center}
\begin{tabular}{|c|c|c|c|c|c|}
\hline
$\alpha$ & Theoretical &  $c(t)$ & $sd$ &  $\hat{c}(t)$ & $\hat{sd}$ \\
  &exponent & $x=\frac{1}{2}$ & $x=\frac{1}{2}$ & $\displaystyle\sup_x$ & $\displaystyle\sup_x$ \\
  \hline
   White noise& 0.5  & 0.6665   &0.0063   &0.6330    & 0.0108   \\
    0.9&0.55 &0.6954 &0.0121 &0.6853 & 0.0130 \\
    0.8&0.6  &0.7548&0.0098 &0.7203 &0.0134 \\
    0.7 &0.65 &0.7512  &0.0089  &0.7508  &0.0186  \\
    0.6& 0.7 &0.8158 &0.0143  &0.8007  &0.0090   \\
    0.5&0.75& 0.8826 &0.0144  &0.8512 &0.0089 \\
    0.4&0.8&0.8987  &0.0100  &0.9112  &0.0113  \\
    0.3&0.85 & 0.9592  & 0.0117 &0.9135  &0.0117  \\
    0.2 & 0.9&0.9891 &0.0116  &0.9563 &0.0147  \\
    0.1&0.95 &1.1797  &0.0114  &1.0219  &0.0120  \\
    \hline
    \end{tabular}
  \end{center}
}
{\small \begin{center}
\begin{tabular}{|c|c|c|c|c|c|}
\hline
$\alpha$ & Theoretical &  $\gamma(x)$& $SD$ &  $\hat{\gamma}(x)$ & $\hat{SD}$ \\
  &exponent & $x=\frac{1}{2}$ & $x=\frac{1}{2}$ & $\displaystyle\sup_x$ & $\displaystyle\sup_x$ \\
 \hline
   White noise& 1.0  & 1.2513  &0.0346   &1.2504  &0.0268    \\
   0.9&1.1 &1.3467 &0.0340 &1.3361  &0.0201  \\
   0.8&1.2  & 1.4347 &0.0336 &1.4251  &0.0211 \\
   0.7 &1.3 & 1.5460 &0.0305  &1.5050  &0.0298  \\
   0.6& 1.4 &1.5869 & 0.0210  &1.5859  &0.0274   \\
   0.5&1.5& 1.6714 &0.0280   &1.6671 &0.0272  \\
   0.4&1.6& 1.7704 &0.0283  &1.7259   &0.0259   \\
   0.3&1.7& 1.8381 &0.0280  &1.7911 &0.0232  \\
   0.2 & 1.8&1.8978 &0.0274  &1.8503 &0.0208  \\
    0.1&1.9 &1.9236 &0.0208  &1.9054 &0.0229  \\
     \hline
     \end{tabular}
\end{center}
}

Finally, since our method applies in the case of non-linear coefficients, we have performed similar computations for
$e(m_i)$, $\hat{e}(m_i)$, $\varepsilon(n_j)$ and $\hat{\varepsilon}(n_j)$ for $1\leq i\leq 13$ and $1\leq j\leq 7$ with
$K=3000$ iterations in the case $\sigma(x)=b(x)=1+0.2\, \cos(x)$. The corresponding results are summarized as follows
{\small \begin{center}
\begin{tabular}{|c|c|c|c|c|c|}
 \hline
 $\alpha$ & Theoretical &  $c(t)$ & $sd$ &  $\hat{c}(t)$ & $\hat{sd}$ \\
  &exponent & $x=\frac{1}{2}$ & $x=\frac{1}{2}$ & $\displaystyle\sup_x$ & $\displaystyle\sup_x$ \\
 \hline
 White noise& 0.5  & 0.4915 &0.0602  &0.5200 &0.0431    \\
  0.8&0.6  & 0.5550  &0.0449 &0.6070  &0.0496 \\
  0.5&0.75&  0.7244 &0.0176  &0.7947 &0.0431 \\
   0.2 & 0.9&0.8607 &0.0225  &0.8571 &0.0429  \\
   \hline
   \end{tabular}
\begin{tabular}{|c|c|c|c|c|c|}
\hline
   $\alpha$ & Theoretical &  $\gamma(x)$& $SD$ &  $\hat{\gamma}(x)$ & $\hat{SD}$ \\
  &exponent & $x=\frac{1}{2}$ & $x=\frac{1}{2}$ & $\displaystyle\sup_x$ & $\displaystyle\sup_x$ \\
    \hline
   White noise& 1.0  & 1.0278   &0.0790  & 0.8263 &0.1056    \\
   0.8&1.2  & 1.3628  &0.0830 & 1.1276  &0.0684 \\
    0.5&1.5& 1.5626  & 0.0710 & 1.5507 &0.0686  \\
    0.2 & 1.8& 1.7351 &0.0708  & 1.4875 &0.0768  \\
   \hline
  \end{tabular}
\end{center}
}
In this semi-linear  case, the speed of convergence is worse and the precision is less than
in the previous linear case.

\noindent {\bf Acknowledgments:} The authors wish to express their
gratitude to Olivier Catoni and Jacques Portes for their helpful
advise when writing the C codes.
\appendix
\section{Appendix.}
\setcounter{equation}{0}
\renewcommand{\theequation}{\thesection.\arabic{equation}}
We start this section with some results concerning the Green
kernel $G_d$ in arbitrary dimension $d\geq 1$. As in the previous
sections, we will suppose that $G_d$ and its discretized versions
are defined with the homogeneous Dirichlet conditions on $\delta
Q$; all the results stated remain  true for the Neumann ones.
\begin{lemme}
\label{A1}
Let $d\geq 1$ and $\alpha\in ]0,2\wedge d[$. There exists some constant $C>0$  depending only on $\alpha$,
such that for all $x$, $x'$ in $Q=[0,1]^d$ and $0< t\leq t'\leq T$~:
\begin{align}
\sup_{y\in Q}\left\| G_d(t,y,\cdot) \right\|_{(\alpha)}^2 \leq& C \, t^{-\frac{\alpha}{2}}\, ,\label{A1(1)}\\
\int_{0}^{+\infty}\left\| G_d(t,x,\cdot)-G_d(t,x',\cdot) \right\|_{(\alpha)}^2dt
\leq&C\, |x-x'|^{2-\alpha}\, ,\label{A1(2)}\\
\sup_{x\in Q} \int_{0}^t\left\| G_d(t'-s,x,\cdot)-G_d(t-s,x,\cdot)\right\|_{(\alpha)}^2ds
\leq&C\, |t'-t|^{1-\frac{\alpha}{2}}\, ,\label{A1(3)}\\
\sup_{x\in Q} \int_{t}^{t'}\left\|G_d(t'-s,x,\cdot)\right\|_{(\alpha)}^2ds
\leq&C\, |t'-t|^{1-\frac{\alpha}{2}}\, .\label{A1(4)}
\end{align}
\end{lemme}
{\bf Proof:} To prove (\ref{A1(1)}), recall the  upper estimate of
$|G_d|$ given in (\ref{A1(0)}).  We remark that  $\ddd \exp\big(-c
|x-y|^2 
t^{-1}\big) \leq \exp\left(-c\, |y-z|^2\,
t^{-1}\right)$
 if $|x-y|\geq |y-z|$, while  $|y-z|^{-\alpha}\geq |x-y|^{-\alpha}$  if $|x-y|\leq |y-z|$.  Hence
(\ref{A1(1)}) follows from
\[ \sup_{x\in Q}\left\|G_d(t,x,\cdot)\right\|_{(\alpha)}^2 
\leq C t^{-d}\left(\int_{0}^{+\infty}\!\!\! e^{-c\frac{u^2}{t}} u^{-\alpha+d-1}\, du\right)
\left(\int_{0}^{+\infty}\!\!\! e^{-c\frac{v^2}{t}} v^{d-1}\, dv\right).\]
We now prove (\ref{A1(2)}) and  set $x'=x+v$. Then, for $0<t\leq |v|^2$, we have
\[
\left\|G_d(t,x,\cdot)-G_d(t,x',\cdot)\right\|_{(\alpha)}^2\leq 2\,
[\, \left\|\,G_d(t,x,\cdot)\,\right\|_\alpha^2 +
\left\|\,G_d(t,x',\cdot)\,\right\|_\alpha^2\, ]\,  .
\]
The change of variables defined by $x-y=|v|\, \eta$, $x-z=|v|\, \xi$ and $t=|v|^2\, s$ in the first integral
(and a similar one with $x'$ instead of $x$ in the second one), combined with  (\ref{A1(1)}),
yields
\begin{align*}
& \int_{0}^{|v|^2}\!\!\left\|G_d(t,x,\cdot)-G_d(t,x',\cdot)\right\|_{(\alpha)}^2 \, dt \leq
 C|v|^{2-\alpha}\, \int_{0}^1\!\! s^{-d}\, ds
\Big\{\dint_{\{|\xi-\eta|\leq|\eta|\}}\!\!e^{-c\frac{|\xi-\eta|^2}{s}} 
 \\
&\quad \times |\xi-\eta|^{-\alpha}e^{-c\frac{|\xi|^2}{s}}d\xi d\eta  
+ \dint_{\{|\xi-\eta|\geq|\eta|\}}\!\! e^{-c\frac{|\eta|^2+ |\xi|^2}{s}}
|\eta|^{-\alpha} d\xi d\eta\Big\}
 \leq 
C |v|^{2-\alpha}. 
\end{align*}
On the other hand, if $t\geq |v|^2$, for every $j\in \{1, \cdots, d\}$,  we use the following well-known
 estimate  (see e.g. \cite{EZ}):
$\left|\frac{\dd}{\dd x_j}G_d(t,x,y)\right|\leq Ct^{-\frac{d+1}{2}}\, \exp\left(-c\frac{|x-y|^2}{t}\right)$  
and the fact that
$|G_d(t,x,y)-G_d(t,x',y)|
\leq \sum_{i=1}^d \left(  \prod_{j=1}^{i-1}|G(t,x_j,y_j)| \right)
 |G(t,x_i,y_i)-G(t,x'_i,y_i)|
\left(  \prod_{j=i+1}^{d}|G(t,x'_j,y_j)| \right)$ for \linebreak[4] $G:=G_1$.
Thus,  Taylor's formula and for every $i\in \{1 \cdots, d\}$ such that $v_i\neq 0$, the change of variables
$x_j-y_j=v_i\, \eta_j$, $x_j-z_j=v_i\, \xi_j$ for $j \leq i$,
$x'_j-y_j=v_i\, \eta_j$, $x'_j-z_j=v_i\, \xi_j$ for $j \geq i+1$, and $t=v_i^2\, s$   yield
\begin{align*}
& \int_{|v|^2}^{+\infty}\left\|G_d(t,x,\cdot)-G_d(t,x',\cdot)\right\|_{(\alpha)}^2 dt
\leq
 C \sum_{i=1}^d 1_{\{v_i\neq 0\}}  |v_i|^{2-\alpha}\int_1^{+\infty} s^{-(d+1)} ds 
\int_{-1}^1 \, d\lambda \\
&\quad \times  
\Big( \prod_{j\neq i}\int_{\mathbb{R}^2} e^{-c\frac{|\eta_j|^2}{s}}
|\eta_j-\xi_j|^{-\alpha_j}
   e^{-c\frac{|\xi_j|^2}{s}} d\xi_j d\eta_j \Big) 
\int_{\mathbb{R}^2} e^{-c\frac{|\eta_i+\lambda|^2}{s}} |\eta_i-\xi_i|^{-\alpha_i}
 e^{-c\frac{|\xi_i+\lambda|^2}{s}} d\xi_i d\eta_j  .
\end{align*}
Splitting again the integrals between $\{|\xi_j-\eta_j|\leq|\eta_j|\}$ and $\{|\xi_j-\eta_j|\geq|\eta_j|\}$
 for $j\neq i$ and $\{|\xi_i-\eta_i|\leq|\eta_i+\lambda|\}$ and $\{|\xi_i-\eta_i|\geq|\eta_i+\lambda |\}$ yields
\[\int_{|v|^2}^{+\infty}\left\|G_d(t,x,\cdot)-G_d(t,x',\cdot)\right\|_{(\alpha)}^2
\leq C\, |v|^{2-\alpha}\int_{1}^{+\infty}s^{-(1+\alpha)}ds\leq C\, |v|^{2-\alpha}.
\] 
This completes the proof of (\ref{A1(2)}).
On the other hand, (\ref{A1(3)}) is obtained using similar arguments and 
the change of variables defined by $t-s=h\, r$, $y-x=\sqrt{h}\, \eta$ and  $z-x=\sqrt{h}\, \xi$ (where $h=t'-t>0$),
Taylor's formula and the estimate
$  \left|\frac{\dd}{\dd t}G_d(t,x,y)\right|
\leq Ct^{-\frac{d+2}{2}}e^{-c\frac{|x-y|^2}{t}}.$ $\Box$

We recall the following well-known set of estimates.  The proofs can be found in \cite{W}
for $d=1$ and are easily deduced for any $d\geq 2$. By convention set $G_d(t,x,y)=0$ if $t\leq 0$.
\begin{lemme}
\label{A2}
For $ x,x'\in Q$, $0\leq t<t'\leq T$ and $\mu\in ]0,1[$,
\begin{eqnarray}
\int_{0}^{+\infty}\hspace{-1mm}\int_Q|G_d(t,x,y)-G_d(t,x',y)|\, dy\, dt&\leq& C\, |x-x'|\label{A2(1)}\, ,\\
\int_{0}^{t'}\hspace{-1mm}\int_Q |G_d(t-s,x,y)-G_d(t'-s,x,y)|\,
 dy\, ds&\leq& C\, |t'-t|^{\mu}\label{A2(3)}\, .
\end{eqnarray}
\end{lemme}

The following technical results are needed to obtain refined estimates for the discretized kernels
$G^n$ and $G^{n,m}$. The proofs, based on simple comparison between series and integrals for
piecewise monotone functions, are omitted.
\begin{lemme}
\label{A4}
For any $c\geq 0$ there exists a constant $C>0$ such that, for $K\geq 0$, $\beta\in [0,1[$, $t> 0$, $a>1$
and $J_{0}\geq 1$,
\begin{align}
& \sum_{j=J_{0}}^{\infty}j^{-\beta}e^{-ctj^2}\leq
C  e^{-ctJ_{0}^2} \left[1+t^{-\frac{1-\beta}{2}}\right] ,
\; \sum_{j=1}^{\infty}j^{K}e^{-ctj^2}\leq
 C \, \left[1+ t^{-\frac{K+1}{2}}\right]\, e^{-ct}\, ,\label{A4(1)}\\
&\sum_{j=J_0}^{\infty} \Big(1+\frac{cTj^2}{m}\Big)^{-1}\leq C m^{\frac{1}{2}} T^{-\frac{1}{2}} ,\; 
\sum_{j=J_0}^{\infty} \Big(1+\frac{cTj^2}{m}\Big)^{-a}\leq 
 \frac{Cm}{J_0 T (a-1)} \left(1+\frac{c T
 J_0^2}{m}\right)^{-a+1} .
\label{A4(3)}
\end{align}
\end{lemme}
The following lemma  bounds the  $\|\; \|_1$ and $\|\;
\|_\alpha$ norms of $(G_d)^n(t,x,.)$.
\begin{lemme}\label{B1n} There exists a constant $C>0$ such that for every $t>0$,
$d\geq 1$, $\lambda >0$ and $0<\alpha<\beta<d\wedge 2$,
\begin{align}
\label{B1(1)}\sup_n\sup_{x\in Q} \|(G_d)^n(t,x,.)\|_1\leq& C\,
\big( 1+ t^{-\lambda}\big) \, e^{-ct}\, ,\\
\label{Gnalpha}\sup_n\sup_{x\in Q}
\|(G_d)^n(t,x,.)\|_{(\alpha)}^2\leq&
C\, t^{-\frac{\beta}{2}}\, e^{-ct}\, ,\\
\label{B1(2)}\sup_{x\in Q} \|(G_d)^n(t,x,.)\|_{(\alpha)}^2\leq&
C\, n^\alpha\, e^{-ct}\, .
\end{align}
\end{lemme}
{\bf Proof:} It suffices to check these inequalities for $x=l/n$, $0\leq l\leq n$. We at first prove them for $d=1$.
Using Abel's summation method, we have, since $ j\longmapsto
\lambda_{j}^n=-4n^2\sin^2\left(\frac{j\pi}{2n}\right)$ is decreasing, for $\kappa=\kappa_n(y)$,
\begin{equation}
|G^n(t,x,y)|
\leq C   e^{\lambda_{1}^nt}
\Big\{\Big|\sin\big(\pi\frac{x-\kappa_n(y)}{2}\big)\Big|^{-1} +
\Big|\sin\big(\pi\frac{x+\kappa_n(y)}{2}\big)\Big|^{-1}\Big\}.
\label{B1(8)}
\end{equation}
Fix $0<\lambda<\frac{1}{2}$ and $t>0$; (\ref{B1(1)}) and
(\ref{A4(1)}) with $\beta=0$, $J_0=1$ yield (\ref{B1(1)}) for
$d=1$.
Let $A_n^i(x)$ be the sets defined by (\ref{B2(7)}) and  for $0\leq l\leq n$, set
$D_{n}^{(i)}(l)= A_n^i(l/n)$.
Then $dx\left(D_{n}^{(i)}(l)\right)\leq\frac{C}{n}$ and, if $y\not\in D_{n}^{(3)}(l)$, one has
$\ddd |x-k_{n}(y)|\geq\frac{2}{3}|x-y|,\;|x+k_{n}(y)|\geq\frac{2}{3}|x+y|$
and similarly, if $y\not\in D_{n}^{(2)}(l)$,
$\ddd |x-k_{n}(y)|\geq\frac{1}{2}|x-y|,\;|x+k_{n}(y)|\geq\frac{1}{2}|x+y|$.
Thus, for every $n\geq 1$ and $0\leq l\leq n$,
$\|G^n(t,l/n,.)\|_{(\alpha)}^2 \leq C(T_1+T_2+T_3)$, where
$T_i$, $1\leq i\leq 3$ is the integral of  $|G^n(t,l/n,y)|\, |y-z|^{-\alpha}\, |G(t,l/n,z)|$ respectively on the
set  $A_1=\{(y,z)\, :\, y\in D_n^{(2)}(l), z\in D_n^{(2)}(l)\} $,
 $A_2=\{(y,z)\, :\, y\in D_n^{(2)}(l)^c,\,  z\in D_n^{(2)}(l)^c,\, |y-z|\leq n^{-1}\}  $
and $A_3=\{(y,z)\, :\, y\in D_n^{(2)}(l)^c,\,  z\in
D_n^{(2)}(l)^c,\, |y-z|\geq n^{-1}\}  $.
Thus, using (\ref{A4(1)})
and (\ref{B1(8)}), we deduce 
upper estimates of $T_i$ for $1\leq i\leq 3$ which imply
(\ref{Gnalpha}) when $d=1$.\\
Again, to prove (\ref{B1(2)}), it suffices to show that 
$\sup_{0\leq l\leq n} \|G^n(t,l/n,.) \|_{(\alpha)}^2 \leq C\,
e^{-ct}\, n^\alpha$. Using the sets $A_i$, $1\leq i\leq 3$, the
inequality (\ref{B1(8)}), the crude estimate
$|G^n(t,x,y)|\leq C \, n\, e^{-ct}$,
 and replacing in products involving two of the terms $|x-y|^{-1}$, $|x-z|^{-1} $
and $|y-z|^{-\frac{\alpha}{2}}$ the largest norm by the smallest one,
a similar computation yields  for $\alpha<\lambda<1$
and $0<\mu<1$,
\[
\|G^n(t,x,.) \|_{(\alpha)}^2 
\leq  C\, \big[\, n^{\lambda+\mu}\, n^{-\lambda-\mu+\alpha}  + n^\alpha\, \big]
\, e^{-ct}\leq C\, n^{\alpha}\, e^{-ct}\, . \]
Since $(G_d)^n(t,x,y)=\prod_{i=1}^d G^n(t,x_i,y_i)$, (\ref{B1(1)})
for $d=1$ immediately yields (\ref{B1(1)}) for any $d$. For $d\geq
2$, and $1\leq i\leq d-1$, set $\alpha_i=\alpha\, 2^{-i}$ and set
$\alpha_d=\alpha_{d-1}$. Then using (\ref{majof}),  the inequality
(\ref{Gnalpha}) (resp. (\ref{B1(2)})) for $d=1$, we deduce 
(\ref{Gnalpha}) (resp. (\ref{B1(2)}))
 for every $d$.
 \qquad$\Box$\\
We now prove a similar result for the norms of
$(G_d)^{n,m}(t,x,.)$.
\begin{lemme}
\label{B1mn}
For every  $\lambda\in]0,\frac{1}{2}[$
and $\beta\in ]\alpha , d\wedge 2[$,
there exist positive  constants $c$ and  $C$ such that for every $t\in ]0,T]$,
\begin{align}\label{B1(3-1)}
\sup_{x\in Q} \| (G_d)^{n,m}(t,x,.)\|_1 \leq & C  
e^{-ct} \big(1+t^{
-\lambda}\big)\\ 
\label{B1(3-2)} \sup_{x\in Q} \| (G_d)^{n,m}(t,x,.)
\|_{(\alpha)}^2 \leq & C e^{-ct}\Big[ (n\wedge \sqrt{m})^\alpha
\wedge\big(1+ t^{-\frac{\beta}{2}} \big) \Big]
+ 
 C   t^{-\frac{\beta}{2}}e^{-ctm} .
\end{align}
\end{lemme}
{\bf Proof : }
For $m\geq 1$, set 
$(\bar{G}_d)^{n,m}(t,x,y):=(G_d)^{n,m}(t,x,y)-(\tilde{G}_d)^{n,m}(t,x,y)$, where 
\begin{equation}
(\tilde{G}_d)^{n,m}(t,x,y)=\!\!\!\!\sum_{\underline{\bf k}\in \{1
,\cdots
 ,  [(n\wedge\sqrt{m})-1)]\}^d} \prod_{i=1}^d \big(1-Tm^{-1}
  \lambda_{k_i}^n\big)^{-\left[\frac{mt}{T}\right]}
\varphi_{\underline{\bf k}}^n(x)  \varphi_{\underline{\bf
k}}(\kappa_{n}(y)) . \label{tildeG}
\end{equation}
Let $(\tilde{G}_1)^{n,m}=\tilde{G}^{n,m}$ and $(\bar{G}_1)^{n,m}=\bar{G}^{n,m}$.
Since $j\rightarrow \big(1-\frac{T}{m}\, \lambda_j^n\big)^{-[\frac{mt}{T}]}$ is decreasing,
$ 
|\tilde{G}^{n,m}(t+T m^{-1},x,y)|\leq C (n\wedge \sqrt{m}) e^{-ct}
$ 
and Abel's summation method yields that  for $x=\frac{l}{n}$ and $\kappa_n(y)=\frac{k}{n}$:
$ 
|\tilde{G}^{m,n}(t+T\, m^{-1},x,y)|\leq
C e^{-ct} \Big[\frac{1}{|x-\kappa_n(y)|}+\frac{1}{x+\kappa_n(y)}+\frac{1}{2-x-\kappa_n(y)}\Big]. $
Finally, since for $j\leq \sqrt{m}$, $\ln\Big( 1+\frac{T}{m}\, 4n^2\, \sin\big( \frac{j\pi}{2n}\big) \Big)
\geq C\, j^2\, T\, m^{-1}$. Using (\ref{A4(1)}) we deduce 
$ 
|\tilde{G}^{n,m}(t+T m^{-1},x,y)|\leq C \sum_{j=1}^{(n\wedge\sqrt{m})-1}e^{-ctj^2}\leq
 C e^{-ct} \big(1+t^{-\frac{1}{2}}\big).
$ 
Thus, repeating the arguments used  to prove (\ref{B1(1)}) -  (\ref{B1(2)})
we deduce that for $\lambda >0$, $0<\alpha <\beta<d\wedge 2$,
\begin{align}\label{B1(4-1)}
\sup_{x\in Q}\| \tilde{G}^{n,m}(t,x,.)\|_1\leq&  C\,e^{-ct}\,  \big( 1+t^{-\lambda}\big)\, ,\\
\sup_{x\in Q}\| \tilde{G}^{n,m}(t,x,.)\|_{(\alpha)}^2\leq& C\,
e^{-ct}\, \Big[ \big(1+t^{-\frac{\beta}{2}}\big) \wedge (n\wedge
\sqrt{m})^\alpha \Big]\, . \label{B1(4-2)}
\end{align}
\noindent We finally give an upper estimate of the norms of
$(\bar{G}_d)^{n,m}(t,x,.)$ and thus we suppose that
$\sqrt{m}<n$. Using (\ref{A4(3)}) we deduce the existence of
positive constants $c,C$ such that for   $t\leq \frac{2T}{m}$,
$\sup_{x,y\in Q} |\bar{G}^{n,m}(t,x,y)|
\leq C \int_{\sqrt{cT2^{-1}}}^{+\infty} (mT)^{-\frac{1}{2}}
 (1+y^2)^{-([\frac{mt}{T}]+1)} dy.$
Hence for $t\leq 2Tm^{-1}$, since $[\frac{mt}{T}]=1$ or $2$, for $x,y\in Q$, $|\bar{G}^{n,m}(t,x,y)|\leq C\,
\leq \sqrt{m}\leq t^{-\frac{1}{2}} $ while for $t\geq 2Tm^{-1}$,
$ |\bar{G}^{n,m}(t,x,y)|\leq C (mT)^{-\frac{1}{2}}
\int_{\sqrt{cT}}^{+\infty} y (1+y^2)^{-([\frac{mt}{T}]+1)} dy \leq C t^{-1} m^{-\frac{1}{2}}
\leq C t^{-\frac{1}{2}} . $
This implies that
\begin{equation}\label{B1(5)}
\sup_{x,y\in Q} |\bar{G}^{n,m}(t,x,y)|\leq C\, \big(1+t^{-\frac{1}{2}}\big)\, .
\end{equation}
Furthermore,  $j \rightarrow (1-T
\lambda^n_j/m)^{-[\frac{mt}{T}]}$ decreases and
$\big(1-T\lambda_{\sqrt{m}}^n/m\big)^{-([\frac{mt}{T}]+1)} \leq
C\,e^{-cTm}$. Hence for $x=l/n$ and $\kappa_n(y)=k/n$  by Abel's
summation method
\begin{equation}\label{B1(6)}
|\bar{G}^{n,m}(t+T\, m^{-1},x,y)|\leq C\, e^{-ctm}\ \Big[\frac{1}{|x-\kappa_n(y)|} + \frac{1}{x+\kappa_n(y)}+
\frac{1}{2-x-\kappa_n(y)}\Big]\, .
\end{equation}
An argument similar to that used to prove (\ref{B1(1)}) implies that for $\lambda \in ]0,\frac{1}{2}[$, there exists
$C>0$ such that for $t\in ]0,T[$,
\begin{equation}\label{B1(6-1)}
\sup_{x\in Q} \| \bar{G}^{n,m}(t,x,.)\|_1\leq C\,
\big(1+t^{-\lambda}\big)\, .
\end{equation}
Finally, for $x=l/n$  let $\bar{D}^i_m(l)=\{ z\in [0,1]: |x-z|\leq i\sqrt{m}, \mbox{\rm or }x+z\leq i\sqrt{m}, \mbox{\rm or }
2-x-z\leq i\sqrt{m}\}$. Then since $n  \geq  \sqrt{m}$, for $y\not\in \bar{D}^3_m(l)$ we deduce that $|x-\kappa_n(y)\geq
\frac{1}{2}|x-y|$, $|x+\kappa_n(y)\geq \frac{1}{2}|x-y|$ and $|2-x-\kappa_n(y)\geq \frac{1}{2}|x-y|$.
Hence, the arguments used to prove (\ref{Gnalpha}) with $d=1$ and (\ref{B1(2)}) show that
$ 
\| \bar{G}^{n,m}(t+T\, m^{-1},x,.) 1_{\bar{D}^3_m(l)^c}\|_{(\alpha)}^2
\leq C e^{-ctm} m^{\frac{\alpha}{2}} .
$ 
Furthermore, (\ref{B1(5)}) and (\ref{B1(6)}) imply that 
 the same upper estimate holds for $ \| \bar{G}^{n,m}(t+T\, m^{-1},x,.) 1_{\bar{D}^3_m(l)}\|_{(\alpha)}^2$.
 These inequalities 
  imply that for $\beta\in ]\alpha,d\wedge 2[$,
\begin{equation}\label{B20}
\sup_{x\in [0,1]} \| \bar{G }^{n,m}(t,x,.)\|_{(\alpha)}^2 \leq C\, e^{-ctm}\, m^{\frac{\alpha}{2}}\,
\big[ 1+(tm)^{-\frac{\beta}{2}}\big] \, .
\end{equation}
Hence (\ref{B20}) and (\ref{B1(4-2)}) imply that  (\ref{B1(3-2)}) holds for $d=1$.
Finally, as in the proof of Lemma \ref{B1n}, (\ref{B1(3-1)}) and (\ref{B1(3-2)}) hold for any $d\geq 1$.
 \qquad$\Box$\\
\noindent We finally prove upper estimates for the norms of time
increments of $(G_d)^n$ and set by convention $(G_d)^n(t,x,.)=0$ if $t\leq 0$.
\begin{lemme}\label{*}
For any 
 $T>0$,  there exists $C>0$ such that for any $h>0$
\begin{align}
&\sup_{n\geq 1}\, \sup_{x\in Q}\, \sup_{t\in [0,T]}
\int_0^{t+h}\|(G_d)^n(t-s,x,.)-(G_d)^n(t+h-s,x,.)\|_1\, ds \leq C\,
h^{\frac{1}{2}}\label{*1}\, ,\\
&\sup_{n\geq 1}\, \sup_{x\in Q}\, \sup_{t\in [0,T]}
\int_0^{t+h}\|(G_d)^n(t-s,x,.)-(G_d)^n(t+h-s,x,.) \|_{(\alpha)}^2\, ds
\leq C\, h^{1-\frac{\alpha}{2}}\label{*4}\, . 
\end{align}
\end{lemme}
{\bf Proof:} Inequality (\ref{B1(1)}) implies that for  $x\in Q$, $\lambda >0$ and $n\geq 1$, $\int_t^{t+h}
  \int_Q |(G_d)^n(t+h-s,x,y)|
 dy ds \leq C \int_0^h s^{-\lambda} ds \leq C h^{1-\lambda}$. 
Using (\ref{B1(2)}) we deduce that
for any $\tilde{c}>0$ and  $h\leq  n^{-2}$,
\begin{equation}\label{*41}
 \int_0^{\tilde{c}h} \|(G_d)^n(s,x,.) \|_{(\alpha)}^2 \, ds \leq \int_0^{\tilde{c}h} n^\alpha \, ds
\leq C\, h^{1-\frac{\alpha}{2}} .
\end{equation}
Suppose that $s\geq n^{-2}$ and  that $d=1$. Then by (\ref{A4(1)}), $|G^n(s,x,y)|\leq C\, (1+s^{-\frac{1}{2}})$. Then  using
(\ref{B1(8)}) and proceeding as in the proof of (\ref{B1(2)}), replacing the sets ${\mathcal A}_n^i(x)$ defined by (\ref{B2(7)}) by the
sets ${\mathcal A}_h^i(x)=\{ y\in [0,1]: |y-x|\leq i\sqrt{h}\; \mbox{\rm or }\; y+x\leq i\sqrt{h}\;
\mbox{\rm or }\; 2-x-y \leq i\sqrt{h}\}$,
since we have assumed that $n^{-1}\leq \sqrt{s}$, we deduce that $\|G^n(s,x,.) \|_{(\alpha)}^2 \leq C\, s^{-\frac{\alpha}{2}}\,
e^{-ct}$. Let $\alpha_k=\alpha\, 2^{-k}$ for $1\leq k\leq d-1$ and $\alpha_d=\alpha_{d-1}$; using the inequality (\ref{majof}),
we deduce that  for $s\geq n^{-2}$, then
$ \|(G_d)^n(s,x,.) \|_{(\alpha )}^2\leq C s^{-\frac{\alpha}{2}} e^{-ct}.$
Hence this inequality and   (\ref{*41})  imply
$ 
\int_t^{t+h}\! \|(G_d)^n(t+h-s,x,.) \|_{(\alpha )}^2 ds \leq C \big[\int_0^{h\wedge n^{-2}} \!\! n^\alpha 
 ds +
\int_{h\wedge n^{-2}}^h s^{-\frac{\alpha}{2}} ds\big] $ $ \leq C h^{1-\frac{\alpha}{2}}$, 
which yields (\ref{*4}).  To complete the proof of  (\ref{*1}) and (\ref{*4}), set $t'=t+h$
and consider the integrals on the interval $[0,t[$. Then for $d=1$ and  
 $x=\frac{l}{ n}$ one has 
$|G^n(t,x,y)-G^n(t',x,y)|= \Big| \sum_{j=1}^{n-1}
 e^{-4tn^2\sin^2\big(\frac{j\pi}{2n}\big)}
 \Big[1-
e^{-4n^2h\sin^2\big(\frac{j\pi}{2n}\big)}
\Big] \varphi_j(x)
\varphi_j(\kappa_n(y)\Big| . $
Thus (\ref{A4(1)}) implies the existence of  $C>0$
 such that for any $n\geq 1$ 
and $x=l\, n^{-1}$,
$ \int_0^t\int_Q|G^n(t-s,x,y)-G^n(t'-s,x,y)| dy ds
\leq C\sum_{j=1}^{n-1}j^{-2}  \big[(j^2 h)\wedge 1\big] \leq
C h^\frac{1}{2},
$ 
which proves (\ref{*1}). The inequality (\ref{*41}) proves that given any $\tilde{c}>0$,
\begin{equation} \label{*5}
\sup_{n\geq 1}\, \sup_{x\in[0,1]}\, \int_0^{\tilde{c}h} \Big( \|G^n(s,x,.)\|_{(\alpha)}^2
+|G^n(s+h,x,.) \|_{(\alpha)}^2 \Big) \, ds
\leq C\, h^{1-\frac{\alpha}{2}}\, .
\end{equation}
Fix $\tilde{c}>0$ large enough, let $t\geq \tilde{c}\, h$ and set
$\Phi(j)=\exp(-4n^2t \sin^2(\frac{j\pi}{2n}))
- \exp(-4n^2t' \sin^2(\frac{j\pi}{2n}))\geq 0$. Then $
\Phi'(j)=2n\sin\big(\frac{j\pi}{n}\big)\big[ t'\, \exp\big(
-4n^2t'\sin^2\big(\frac{j\pi}{2n}\big)\big) -t \, \exp\big(
-4n^2t'\sin^2\big(\frac{j\pi}{2n}\big)\big) \big] $. Then the
arguments used to estimate $\Phi_2(x)$ and then $|T_2(t,x,y)|$ in
the proof of Lemma \ref{B3} show that there exists  $C>0$ such
that for any $s\in [\tilde{c}h, T]$, $x=l\, n^{-1}$ and $y\in
[0,1]$,
$ 
 \big| \sum_{j=[h^{-\frac{1}{2}}]}^{n-1} \big( e^{\lambda_j^n s }- e^{\lambda_j^n (s+h)}\big)
\varphi_j(x)\, \varphi_j(\kappa_n(y))\big|
\leq  C\, \exp\big(\lambda_{[h^{-\frac{1}{2}}]}^n s\big)
\big[ \frac{1}{|x-\kappa_n(y)|} +
\frac{1}{x+\kappa_n(y)}
+\frac{1}{2-x-\kappa_n(y)}\big]\, , 
$  
while (\ref{A4(1)}) implies that
 $ 
\big| \sum_{j=[h^{-\frac{1}{2}}]}^{n-1} \big( e^{\lambda_j^n s }-
 e^{\lambda_j^n (s+h)}\big)\varphi_j(x) \varphi_j(\kappa_n(y)\big|
\leq C s^{-\frac{1}{2}} 
\exp\big(\lambda_{[h^{-\frac{1}{2}}]}^n
s\big) $. 
Using again the sets ${\mathcal A}^i_h(x)$, we deduce
that for any $s\in [\tilde{c}h,T]$, 
$ 
\sup_{n\geq 1} \sup_{x\in [0,1]} \|G^n(s,x,.)-G^n(s+h,x,.) \|_{(\alpha)}^2 \leq C  h^{1-\frac{\alpha}{2}}.
$ 
Thus,
$ \sup_{n\geq 1} \sup_{x\in [0,1]} \int_{\tilde{c}h}^T \|G^n(s,x,.)-G^n(s+h,x,.) \|_{(\alpha )}^2 ds \leq C
 h^{1-\frac{\alpha}{2}}.$
This inequality and  (\ref{*5}) 
yield (\ref{*4}) for $d=1$. 
We extend the lemma in any dimension $d\geq 1$ as in the proof of Lemma \ref{B1n}. \qquad$\Box$ 


\end{document}